\documentclass[english]{amsart}
\pdfoutput=1
\usepackage[utf8]{inputenc}
\usepackage[T1]{fontenc}
\usepackage{babel}

\usepackage{amssymb}
\usepackage{rotating}
\usepackage[all, cmtip]{xy}

\usepackage{amsthm}
\newtheorem{theorem}{Theorem}[section]
\newtheorem{corollary}[theorem]{Corollary}
\newtheorem{lemma}[theorem]{Lemma}
\newtheorem{proposition}[theorem]{Proposition}

\theoremstyle{definition}
\newtheorem{definition}[theorem]{Definition}
\newcommand{\keyterm}[1]{\textbf{#1}}
\newtheorem{example}[theorem]{Example}

\theoremstyle{remark}
\newtheorem{remark}[theorem]{Remark}
\numberwithin{equation}{section}

\usepackage{paralist}
\setdefaultenum{\rm 1)}{\rm a)}{\rm i)}{\rm A)}

\usepackage{varioref}
\labelformat{chapter}{Chapter~#1}
\labelformat{figure}{Figure~#1}
\labelformat{section}{Section~#1}
\labelformat{subsection}{Section~#1}
\labelformat{table}{Table~#1}

\usepackage[dvipsnames]{xcolor}
\usepackage{sagetex}
\lstdefinestyle{SageTEXAdjustment}{numbers=none, xleftmargin=1ex, xrightmargin=1ex, framesep=1ex, frame=rl, breakatwhitespace=True, basicstyle={\ttfamily\small}, backgroundcolor=\color{Goldenrod!25}, framerule=0pt}
\lstdefinestyle{SageInput}{style=DefaultSageInput, style=SageTEXAdjustment}
\lstdefinestyle{SageOutput}{style=DefaultSageOutput, style=SageTEXAdjustment, basicstyle={\ttfamily\small\bfseries}}
\newcommand{\sagecommand}[1]{\lstinline!#1!}

\usepackage{tikz}
\usepackage{tkz-graph}
\usepackage{tkz-berge}
\usetikzlibrary{arrows,shapes}

\usepackage[breaklinks=true, citecolor=dbluecolor, colorlinks=true, hyperfootnotes=false, linkcolor=dbluecolor, urlcolor=dbluecolor]{hyperref}

\newcommand{\mpart}{\mathrm{part}}

\renewcommand{\epsilon}{\varepsilon}
\renewcommand{\kappa}{\varkappa}
\renewcommand{\phi}{\varphi}
\renewcommand{\rho}{\varrho}

\newcommand{\pair}[1]{\left\langle #1 \right\rangle}
\newcommand{\set}[1]{\left\{#1\right\}}

\newcommand{\tld}[1]{\widetilde{#1}}

\newcommand{\CC}{\ensuremath{\mathbb{C}}}
\newcommand{\PP}{\ensuremath{\mathbb{P}}}
\newcommand{\QQ}{\ensuremath{\mathbb{Q}}}
\newcommand{\RR}{\ensuremath{\mathbb{R}}}
\newcommand{\WP}{\ensuremath{\mathbb{WP}}}
\newcommand{\ZZ}{\ensuremath{\mathbb{Z}}}

\newcommand{\cF}{\ensuremath{\mathcal{F}}}
\newcommand{\cO}{\ensuremath{\mathcal{O}}}

\newcommand{\cE}{\ensuremath{\mathcal{E}}}
\newcommand{\cA}{\ensuremath{\mathcal{A}}}

\DeclareMathOperator{\Hom}{Hom}
\DeclareMathOperator{\NS}{NS}

\DeclareMathOperator{\Spec}{Spec}

\begin{document}
\begin{sagesilent}
toric_plotter.options(font_size=10, ray_thickness=1, point_size=7, ray_label="\\rho")
sage.geometry.point_collection.PointCollection.output_format("separated column matrix")
\end{sagesilent}

\title{The 14th case VHS via K3 fibrations}

\author[A. Clingher]{Adrian Clingher}
\address{
Department of Mathematics and Computer Science, University of Missouri -- St. Louis, St. Louis, MO, 63121, USA}
\email{clinghera@umsl.edu}
\thanks{A. Clingher was supported by Simons Foundation grant no. 208258 and by a Bitdefender Invited Professor scholarship from IMAR}

\author[C.F. Doran]{Charles F. Doran}
\address{
Department of Mathematical and Statistical Sciences, 632 CAB, University of Alberta, Edmonton, Alberta, T6G 2G1, Canada}
\email{charles.doran@ualberta.ca}
\thanks{C. F. Doran and A. Y. Novoseltsev were supported by the Natural Sciences and Engineering Resource Council of Canada (NSERC), the Pacific Institute for the Mathematical Sciences, and the McCalla Professorship at the University of Alberta}

\author[J. Lewis]{Jacob Lewis}
\address{Fakult\"{a}t f\"{u}r Mathematik, Universit\"{a}t Wien, Garnisongasse 3/14, 1090 Wien, Austria}
\email{Jacobml@u.washington.edu}
\thanks{J. Lewis was supported in part by NSF grant  OISE-0965183}

\author[A.Y. Novoseltsev]{Andrey Y. Novoseltsev}
\address{
Department of Mathematical and Statistical Sciences, 632 CAB, University of Alberta, Edmonton, Alberta, T6G 2G1, Canada}
\email{novoselt@ualberta.ca}

\author[A. Thompson]{Alan Thompson}
\address{
Fields Institute, 222 College Street, Toronto, Ontario, M5T 3J1, Canada}
\email{amthomps@ualberta.ca}
\thanks{A. Thompson was supported in part by NSERC and in part by a Fields Institute Ontario Postdoctoral Fellowship with funding provided by NSERC and the Ontario Ministry of Training, Colleges
and Universities}

\date{\today}

\begin{abstract}
We present a study of certain singular one-parameter subfamilies of Calabi-Yau threefolds realized as anticanonical hypersurfaces or complete intersections in toric varieties. Our attention to these families is motivated by the Doran-Morgan classification of variations of Hodge structure which can underlie families of Calabi-Yau threefolds with $h^{2,1} = 1$ over the thrice-punctured sphere. We explore their torically induced fibrations by $M$-polarized K3 surfaces and use these fibrations to construct an explicit geometric transition between an anticanonical hypersurface and a nef complete intersection through a singular subfamily of hypersurfaces. Moreover, we show that another singular subfamily provides a geometric realization of the missing ``14th case'' variation of Hodge structure from the Doran-Morgan list.
\end{abstract}

\maketitle

\tableofcontents

\section{Introduction}

In their paper \cite{DoranMorgan:2006}, Doran and Morgan give a classification of the possible variations of Hodge structure that can underlie families of Calabi-Yau threefolds with $h^{2,1} = 1$ over the thrice-punctured sphere. They find fourteen possibilities. At the time of publication of \cite{DoranMorgan:2006}, explicit families of Calabi-Yau threefolds realising thirteen of these cases were known and are given in \cite[Table 1]{DoranMorgan:2006}. The aim of this paper is to give a geometric example which realizes the fourteenth and final case (henceforth known as the \keyterm{14th case}) from their classification, and to study its properties.

By analogy with other examples (see \cite[Section 4.2]{DoranMorgan:2006}), one might expect that the 14th case variation of Hodge structure should be realized by the mirror of a complete intersection of bidegree $(2,12)$ in the weighted projective space $\WP(1,1,1,1,4,6)$. However, this ambient space is not Fano, so the Batyrev-Borisov mirror construction cannot be applied to obtain such a mirror family.

Instead, Kreuzer and Sheidegger \cite{KKRS:2005} suggest working with a slightly different ambient space, given by a non-crepant blow up of $\WP(1,1,1,1,4,6)$. However, the complete intersection Calabi-Yau threefold of bidegree $(2,12)$ in this ambient space has $h^{1,1} = 3$, so its mirror will have $h^{2,1} = 3$, making it unsuitable as a candidate for the 14th case on the Doran-Morgan list.

We take as the starting point for our discussion this mirror family with $h^{2,1} = 3$, which we denote by $Y$. We begin by constructing it explicitly and exploring its properties, noting that it has a singular one-parameter subfamily $Y_1$ with GKZ series matching that predicted for the 14th case. Moreover, we find that $Y$ has a torically induced fibration by K3 surfaces given as anticanonical hypersurfaces in the three-dimensional space polar to $\WP(1,1,4,6)$.

The existence of this fibration suggests that a model for the 14th case variation of Hodge structure could also be constructed by considering Calabi-Yau threefolds admitting fibrations by anticanonical hypersurfaces in the space polar to $\WP(1,1,4,6)$. A computer search for such fibrations yields a second candidate family, $Z$, which is mirror to a hypersurface of degree $24$ in the weighted projective space $\WP(1,1,2,8,12)$. The threefolds in this family also have $h^{2,1} = 3$.

This poses the natural question of whether the two families $Y$ and $Z$ are related. A careful study of the geometry of each shows that the generic fibre of their K3 fibrations is polarized by the lattice $M := H \oplus E_8 \oplus E_8$. Lattice polarized K3 surfaces of this type were studied by \cite{ClingherDoran:2007} and \cite{CDLW:2009}; using their results we are able to show that the K3 fibrations on $Y$ and $Z$ are very closely related. With this relationship as a guide, we proceed to find an explicit geometric transition between the two, by degenerating $Z$ to a singular subfamily $Z_2$ and blowing up to obtain $Y$.

Having shown that the families $Y$ and $Z$ are in fact two halves of the same picture, we turn our focus back to the family $Y$ and, in particular, its singular subfamily $Y_1$. The main results of this paper are Lemma \ref{lemma:Y1MHS} and Theorem \ref{thm:14thcase}, which describe the mixed Hodge structure on the cohomology groups of $Y_1$. Together they imply:

\begin{theorem} The mixed Hodge structure on the cohomology groups $H^i(Y_1)$ is pure unless $i=3$ and the weight filtration $W_{\bullet}$ on $H^3(Y_1)$ has the following description:
\begin{align*} \mathrm{Gr}^{W_{\bullet}}_kH^3(Y_1) &= 0 \quad \mathrm{for}\ k \neq 2,3\\
\dim_{\CC}W_2(H^3(Y_1)) & = 2
\end{align*}
Furthermore, as $Y_1$ varies in its one-parameter subfamily, the third graded piece $\mathrm{Gr}^{W_{\bullet}}_3H^3(Y_1)$ admits a pure variation of Hodge structure of weight $3$ and type $(1,1,1,1)$, which realizes the 14th case variation of Hodge structure.
\end{theorem}

A crucial step in the proof of this theorem involves quotienting the threefold $Y_1$ by a certain involution and resolving to obtain a new family of Calabi-Yau threefolds $W$. This new family is closely related to $Y_1$ and is interesting in its own right. The K3 fibration on $Y_1$ induces a fibration of $W$ by Kummer surfaces, via the mechanism discussed in \cite{ClingherDoran:2007}. We digress to show how $W$ may be constructed directly as a family of Kummer surfaces using the techniques of \cite{DHNT:2013}; this provides a great deal of insight into the geometry of this new threefold.

Finally, we conclude by discussing the mirrors of the families $Y$ and $Z$ considered here, and compute a mirror for the geometric transition between them. In particular, we find that the singular subfamily $Y_1$ may naturally be seen as the mirror of a complete intersection of bidegree $(2,12)$ in $\WP(1,1,1,1,4,6)$, as originally expected. We note, however, that $Y_1$ does not admit a Calabi-Yau resolution, despite the fact that its mirror admits a Calabi-Yau smoothing. This provides a counterexample to a conjecture of Morrison \cite{Morrison:1999}.
\medskip

The structure of this paper is as follows. In \ref{sect:toric}, we give an overview of relevant results from toric geometry that will be used throughout the rest of the paper; this also serves the function of setting up the notation that will be used in what follows. Then in \ref{sect:models} we construct the families $Y$ and $Z$ and study some of their basic properties, these are summarized by Propositions \ref{prop:Yproperties} and \ref{prop:Zproperties}. In particular, Proposition \ref{prop:Yproperties} shows that $Y$ has a subfamily $Y_1$ with GKZ series matching that predicted for the 14th case variation of Hodge structure.

We begin \ref{sect:matching} with a detailed study of the K3 fibrations on $Y$ and $Z$, using results of \cite{ClingherDoran:2007} and \cite{CDLW:2009}. From this, we show that the parameters describing the K3 fibration on $Y$ match with those describing the K3 fibration on a certain subfamily $Z_2$ of $Z$, and those describing the K3 fibration on the subfamily $Y_1$ match with those describing the K3 fibration on a further subfamily $Z_1$ of $Z_2$. These subfamilies $Y_1$, $Z_1$ and $Z_2$ of $Y$ and $Z$ are studied in \ref{sect:singular}; their properties are summarized by Propositions \ref{prop:Z1Z2properties} and \ref{prop:Y1nodal}. Using these results, we show in Proposition \ref{prop:geomtrans} that the families $Y$ and $Z$ are in fact related by a geometric transition, through the subfamily $Z_2$, and that this transition also relates the subfamilies $Y_1$ and $Z_1$.

\ref{section:involutions} begins with a result, Lemma \ref{lemma:Y1MHS}, which describes the mixed Hodge structure on the cohomology groups of the subfamily $Y_1$. However, there is an ambiguity in the description of the weight filtration on $H^3(Y_1)$ given by this lemma, which prevents us from drawing any conclusions. To resolve this, we study the action of certain involutions on $Y$, which restrict to the subfamily $Y_1$. The resolved quotient of $Y_1$ by one such involution gives a new Calabi-Yau threefold, $W$, with geometry closely related to that of $Y_1$. Proposition \ref{prop:WMHS} and Corollary \ref{corollary:Wsympres} describe this geometry; from them we are able to deduce Theorem  \ref{thm:14thcase} and Corollary \ref{corollary:Y1sympres}, which complete the description of the subfamily $Y_1$. In particular, we find that $Y_1$ realizes the 14th case variation of Hodge structure that we have been seeking.

In \ref{section:forward construction} we digress to discuss the geometry of the Calabi-Yau threefolds $W$. The K3 fibration on $Y_1$ induces a K3 fibration on $W$, the general fibre of which is a Kummer surface. We use this K3 fibration, along with the methods of \cite{DHNT:2013}, to construct a model for $W$ directly. The geometric nature of this construction provides a great deal of insight into the geometry of $W$ and demonstrates the utility of the methods of \cite{DHNT:2013} in a practical example.

Finally, in \ref{section:mirrormoduli} we compute mirrors for $Y$ and $Z$, their subfamilies $Y_1$, $Z_1$ and $Z_2$, and exhibit a mirror for the geometric transition between them. In particular, we find that $Y_1$ may be seen as a mirror to a complete intersection of degree $(2,12)$ in $\WP(1,1,1,1,4,6)$, as expected. Finally, we show that $Y_1$ and its mirror provide a counterexample to a conjecture of Morrison \cite{Morrison:1999}.
\medskip

Most of the computer-aided computations performed for this work were done in the free open source mathematics software Sage~\cite{Sage}, specifically using the toric geometry framework developed by Braun and Novoseltsev~\cite{BraunNovoseltsev:toric_variety}. At the time of this writing, it is necessary to use a small extra patch\footnote{Available as a git branch at \url{https://github.com/novoselt/sage/commits/Cayley}}, but we plan to integrate it into the official release of Sage. An interested reader may easily reproduce all computations used in this article following the presented code, as well as adjust them to suit his or her own needs. Typesetting of code snippets and some of the complicated expressions was greatly facilitated by Sage\TeX~\cite{SageTeX}.
\medskip

\subsection{Acknowledgements} The authors would like to thank Andrew Harder and Ursula Whitcher for their assistance with some of the calculations in \ref{section:forward construction}. A. Clingher would like to thank the Institute of Mathematics ``Simion Stoilow'' of the Romanian Academy in Bucharest, the Institut des Hautes \'{E}tudes Scientifique in Bures-sur-Yvette and the Max Planck Institute for Mathematics in Bonn for their hospitality and support during the academic year 2012-2013. A portion of this work was completed while A. Thompson was in residence at the Fields Institute Thematic Program on Calabi-Yau Varieties: Arithmetic, Geometry and Physics; he would like to thank the Fields Institute for their support and hospitality.

\section{Toric Geometry Review and Notation}\label{sect:toric}

In this section we recall the basic facts of toric geometry in order to set up notation that will be used in the rest of the paper; we refer to the wonderful book~\cite{CoxLittleSchenck:2011} for details. We also briefly describe Batyrev's construction of mirror families of Calabi-Yau anticanonical hypersurfaces in toric varieties~\cite{Batyrev:1994} and its generalization to nef complete intersections by Batyrev and Borisov~\cite{BatyrevBorisov:1996CYCIinTV}. We conclude with a discussion of fibrations of Calabi-Yau varieties that can be obtained from natural fibrations of their ambient spaces.

\subsection{Toric Varieties}

Let $M$ and $N$ be dual lattices of rank $n$, i.e. $M \simeq N \simeq \ZZ^n$ as free Abelian groups, where $M$ is identified with $\Hom_\ZZ(N, \ZZ)$ and $N$ with $\Hom_\ZZ(M, \ZZ)$ using the natural pairing $\pair{\cdot, \cdot} \colon M \times N \to \ZZ$. The real vector space associated to a lattice will be denoted by the subscript $\RR$, e.g. the real vector space associated to the lattice $N$ is $N_\RR = N \otimes_\ZZ \RR.$ The torus of a lattice $N$ is defined to be $T_N = N \otimes_\ZZ \CC^* \simeq (\CC^*)^n$, it is an affine variety with group structure corresponding to componentwise multiplication in $(\CC^*)^n$ and character lattice $M$.

A cone $\sigma \subset N_\RR$ will always be convex, polyhedral, and rational. Its dual cone will be denoted by $\sigma^\vee \subset M_\RR$. The affine toric variety corresponding to a cone $\sigma$ will be denoted by $U_\sigma = \Spec\left(\CC[\sigma^\vee \cap M]\right)$.

For a fan $\Sigma$ in $N_\RR$, we will denote by $X_\Sigma$ the corresponding toric variety and by $\Sigma(k)$ the set of its $k$-dimensional cones. In particular, $\Sigma(1)$ is the set of rays and for each ray $\rho \in \Sigma(1)$ we have a primitive integral generator $v_\rho \in N$, a homogeneous coordinate $z_\rho$, and a prime torus-invariant Weil divisor $D_\rho$ given by $\set{z_\rho = 0}$.

Any toric divisor $D = \sum_\rho a_\rho D_\rho$ has an associated polyhedron given by
\begin{gather*}
P_D = \set{u \in M_\RR : \pair{u, v_\rho} \geqslant - a_\rho \text{ for all } \rho \in \Sigma(1)}.
\end{gather*}
For any $m \in P_D \cap M$, the $D$-homogenization of the character $\chi^m$ is the monomial
\begin{gather*}
z^{\pair{m, D}} = \prod_\rho z_\rho^{\pair{m, v_\rho} + a_\rho},
\end{gather*}
and the global sections of $D$ are
\begin{gather*}
\Gamma(X_\Sigma, \cO_{X_\Sigma}(D)) = \bigoplus_{m \in P_D \cap M} \CC \chi^{m}.
\end{gather*}

Of particular interest to us are the toric canonical divisor $K_{X_\Sigma} = - \sum_\rho D_\rho$ and the global sections of the anticanonical divisor, in the following special case.

\subsection{Gorenstein Fano Toric Varieties}

Gorenstein Fano toric varieties are in bijection (after taking into account equivalences) with reflexive polytopes. A reflexive polytope is a full-dimensional lattice polytope $\Delta \subset M_\RR$ that contains the origin in its interior and has polar
\begin{gather*}
\Delta^\circ = \set{v \in N_\RR : \pair{u, v} \geqslant - 1 \text{ for all } u \in \Delta} \subset N_\RR
\end{gather*}
that is also a lattice polytope. Since $\left(\Delta^\circ\right)^\circ = \Delta$, there is a polar duality between reflexive polytopes in $M_\RR$ and reflexive polytopes in $N_\RR$, and even an inclusion-reversing bijective correspondence between the faces of $\Delta$ and $\Delta^\circ$. 

Given a reflexive polytope $\Delta \subset M_\RR$, its normal fan in $N_\RR$ will be denoted by $\Sigma_\Delta$. To avoid towers of subscripts, we will simply write $X_\Delta$ instead of $X_{\Sigma_\Delta}$ for the toric variety associated to this fan. The vertices of $\Delta^\circ$ coincide with the primitive integral generators of rays of $\Sigma_\Delta$, so one can also think of $\Sigma_\Delta$ as the face fan of $\Delta^\circ$ (with cones over faces of $\Delta^\circ$).

\begin{definition}
Let $\Delta \subset M_\RR$ be a reflexive polytope and let $\Sigma$ be a subdivision of $\Sigma_\Delta$. If all rays of $\Sigma$ are generated by (some of) the boundary lattice points of $\Delta^\circ$, it is a \keyterm{crepant subdivision}. If $\Sigma$ is also simplicial and the corresponding toric variety $X_\Sigma$ is projective, it is a \keyterm{projective crepant subdivision}. If furthermore \emph{all} boundary lattice points of $\Delta^\circ$ generate rays of $\Sigma$, it is a \keyterm{maximal projective crepant subdivision}.
\end{definition}

Such subdivisions correspond to \emph{maximal projective crepant partial desingularizations (MPCP-desingularizations)}, introduced by Batyrev in~\cite{Batyrev:1994}. They lead to a very concrete construction of Calabi-Yau orbifolds.

\begin{theorem}
Let $\Delta$ be a reflexive polytope of dimension $n$. A generic anticanonical hypersurface in $X \!=\! X_\Delta$, i.e. a generic section $f \in \Gamma(X, \cO_X(-K_X))$, is a Calabi-Yau variety of dimension $n-1$. A generic anticanonical hypersurface in $X = X_\Sigma$, where $\Sigma$ is a projective crepant subdivision of $\Sigma_\Delta$, is a Calabi-Yau orbifold of dimension $n-1$.
\end{theorem}

\begin{proof}
See \cite[Proposition~4.1.3]{CoxKatz:1999}.
\end{proof}

For a reflexive polytope $\Delta \subset M_\RR$, projective crepant subdivisions of $\Sigma_\Delta$ correspond to special triangulations of the boundary of $\Delta^\circ$. For any such subdivision the polytope of the anticanonical divisor is easily seen to be $\Delta$ itself. This means that in the equations of the anticanonical Calabi-Yau hypersurfaces, the variables correspond to lattice points of $\Delta^\circ$ while the (coefficients of) monomials correspond to lattice points of $\Delta$:
\begin{gather*}
\sum_{m \in \Delta \cap M}
a_m
\prod_{v_\rho \in \partial \Delta^\circ \cap N}
z_\rho^{\pair{m, v_\rho} + 1}
=
0.
\end{gather*}

Polar duality of reflexive polytopes means that the roles of $\Delta$ and $\Delta^\circ$ in the above discussion can be reversed, leading to another family of hypersurfaces in another toric variety. Batyrev showed~\cite{Batyrev:1994} that if $\dim \Delta = 4$, then generic anticanonical hypersurfaces in MPCP-desingularizations of $X_\Delta$ are \emph{smooth} Calabi-Yau threefolds and the exchange $\Delta \leftrightsquigarrow \Delta^\circ$ corresponds to the exchange $h^{1,1} \leftrightsquigarrow h^{2,1}$ of the Hodge numbers of the anticanonical hypersurfaces of the two families, making them candidates for mirror pairs.

\subsection{Nef Complete Intersections}\label{sec:NefCI}

Batyrev and Borisov generalized the  construction of Calabi-Yau varieties as anticanonical hypersurfaces in toric varieties to the case of complete intersections associated to nef-partitions of reflexive polytopes~\cite{Borisov:1993, BatyrevBorisov:1996CYCIinTV}.

\begin{definition}
Let $\Delta \subset M_\RR$ be a reflexive polytope. A \keyterm{nef-partition} is a decomposition of the vertex set $V$ of $\Delta^\circ \subset N_\RR$ into a disjoint union
\begin{align*}
V = V_0 \sqcup V_1 \sqcup \dots \sqcup V_{r-1}
\end{align*}
such that all divisors $E_i = \sum_{v_\rho \in V_i} D_\rho$ are Cartier. Equivalently, let $\nabla_i \subset N_\RR$ be the convex hull of the vertices from $V_i$ and the origin. These polytopes form a nef-partition if their Minkowski sum $\nabla \subset N_\RR$ is a reflexive polytope. 

The \keyterm{dual nef-partition} is formed by the polytopes $\Delta_i \subset M_\RR$ of the $E_i$, which give a decomposition of the vertex set of $\nabla^\circ \subset M_\RR$; their Minkowski sum is $\Delta$.
\end{definition}

``Nef-partition'' may refer to any of the following decompositions:
\begin{compactenum}
\item $V(\Delta^\circ)$ into a disjoint union of $V_i$,
\item $\Delta$ into a Minkowski sum of $\Delta_i$,
\item the anticanonical divisor of $X_\Delta$ into a sum of $E_i$.
\end{compactenum}
Each of these decompositions can be easily translated into the others. Some care should be taken only to avoid mixing a nef-partition and its dual.

It follows from the definition that polar duality of reflexive polytopes switches convex hull and Minkowski sum for dual nef-partitions:
\begin{align*}
\allowdisplaybreaks
\Delta^\circ
&=
\mathrm{Conv} \left(\nabla_0, \nabla_1, \dots, \nabla_{r-1}\right), \\
\nabla^{\phantom{\circ}}
&=
\nabla_0 + \nabla_1 + \dots + \nabla_{r-1}, \displaybreak[0]\\[1ex]
\Delta^{\phantom{\circ}}
&=
\Delta_0 + \Delta_1 + \dots + \Delta_{r-1}, \\
\nabla^\circ
&=
\mathrm{Conv} \left(\Delta_0, \Delta_1, \dots, \Delta_{r-1}\right).
\end{align*}

Given a nef-partition of an $n$-dimensional reflexive polytope consisting of $r$-parts, generic sections of the divisors $E_i$ determine an $(n-r)$-dimensional complete intersection Calabi-Yau variety. In~\cite{BatyrevBorisov:1996STHN} Batyrev and Borisov show that such varieties corresponding to dual nef-partitions have mirror-symmetric \emph{stringy} Hodge numbers, which we will denote by $h_{st}^{p,q}(Y)$.

\subsection{Torically Induced Fibrations}
\label{subsection:torically induced fibrations}

Let $N$ and $N'$ be lattices, $\Sigma$ be a fan in $N_\RR$, and $\Sigma'$ be a fan in $N'_\RR$. There is a bijection between toric morphisms $\phi\colon X_\Sigma \to X_{\Sigma'}$ (i.e. morphisms of toric varieties preserving the group structure of their tori) and fan morphisms $\tld{\phi}\colon \Sigma \to \Sigma'$ (i.e. lattice homomorphisms $\tld{\phi}\colon N \to N'$ compatible with the fan structure: the linear extension $\tld{\phi}_\RR \colon N_\RR \to N'_\RR$ of $\tld{\phi}$ maps each cone $\sigma \in \Sigma$ into a single cone $\sigma' \in \Sigma'$).

Consider the special case when the lattice homomorphism is surjective, i.e. we have an exact sequence of lattices
\begin{gather*}
0 \to N_0 \to N \xrightarrow{\tld{\phi}} N' \to 0,
\end{gather*}
where $N_0 = \ker \tld{\phi}$. Let $\Sigma_0 = \set{\sigma \in \Sigma : \sigma \subset (N_0)_\RR}$. We can consider $\Sigma_0$ either as a fan in $N_\RR$ or as a fan in $(N_0)_\RR$, giving two corresponding toric varieties $X_{\Sigma_0, N}$ (a dense subset of $X_\Sigma$) and $X_{\Sigma_0, N_0}$.

As discussed in \cite[\S~3.3]{CoxLittleSchenck:2011}, there is a clear relation between these two varieties,
\begin{gather*}
X_{\Sigma_0, N} \simeq X_{\Sigma_0, N_0} \times T_{N'}.
\end{gather*}
In fact $X_{\Sigma_0, N} = \phi^{-1}(T_{N'})$, so part of $X_\Sigma$ is a fibre bundle over $T_{N'}$, with fibres $X_{\Sigma_0, N_0}$. Moreover, if $\Sigma$ is \emph{split by $\Sigma'$ and $\Sigma_0$}, then the whole of $X_\Sigma$ is a fibre bundle over $X_{\Sigma'}$ (see \cite[Definition~3.3.18 and Theorem~3.3.19]{CoxLittleSchenck:2011}). However, splitting is a very strong condition on fans, so instead of imposing it we will work with more general fibrations than fibre bundles.

\begin{definition}
Let $\phi\colon X \to Y$ be a morphism between two varieties. Then $\phi$ is a \keyterm{fibration} if it is surjective and all of its fibres have the same dimension $\dim X - \dim Y$.
\end{definition}

In \cite[Proposition~2.1.4]{HuLiuYau:2002} the authors provide a detailed description of the fibres of \emph{arbitrary} toric morphisms\footnote{Note that while the notation in~\cite{HuLiuYau:2002} is very similar to ours, the authors sometimes \emph{implicitly} assume that toric varieties in question are  \emph{complete}.}. However, we are primarily interested in fibrations ${\phi\colon X_\Sigma \to X_{\Sigma'}}$, as they may induce fibrations $\phi\big|_Y \colon Y \to X_{\Sigma}$ of Calabi-Yau subvarieties $Y$ realized as anticanonical hypersurfaces or nef complete intersections in $X_\Sigma$. This fibration condition may prevent $\Sigma$ from being ``too refined'', leading to singularities of $Y$, but it may still be possible to compose $\phi$ with a crepant resolution of singularities in such a way that $Y$ becomes smooth and the restriction of the composition to $Y$ is still a fibration. 

We now describe a strategy for searching for such toric fibrations. Let $\Delta \subset M_\RR$ be a reflexive polytope and let $\Sigma$ be a crepant subdivision of $\Sigma_\Delta$. Suppose that $\tld{\phi}\colon \Sigma \to \Sigma'$ is a fibration (meaning that $\phi\colon X_\Sigma \to X_{\Sigma'}$ is a fibration in the above sense). As before, its fibre is determined by the subfan $\Sigma_0$ of $\Sigma$ in the sublattice $N_0 = \ker \tld{\phi}$ of $N$. Since we would like the fibres of Calabi-Yau subvarieties to be lower-dimensional Calabi-Yau varieties, it is natural to require that $\Sigma_0$ \emph{is also associated to a reflexive polytope}, i.e. that it is a crepant subdivision of $\Sigma_\nabla$ for some $\nabla \subset (M_0)_\RR$, where $M_0$ is the dual lattice of $N_0$. In this case $\nabla^\circ$ is a ``slice'' of $\Delta^\circ$ by a linear subspace, so one can search for such slices of $\Delta^\circ$ and take $\tld{\phi}$ to be the projection along the linear subspace of $\nabla^\circ$. 

We can also reformulate this problem in dual terms: $M_0 = M / (N_0)^\perp$ and the condition that $\nabla^\circ$ is inside $\Delta^\circ$ implies that the image of $\Delta$ in $M_0$ is inside $\nabla$. So alternatively one can look for ``projections'' of $\Delta$ that are reflexive. Here we recall that the origin is the only interior lattice point of any reflexive polytope, so all lattice points of $\Delta$ must be projected onto either the origin or the boundary of the projection. Due to this restriction a ``large'' $\Delta$ with many lattice points is less likely to have any fibrations than a ``small'' one.

\section{Models} \label{sect:models}

\subsection{Complete Intersections}
\label{subsection:CI model} 

In~\cite{DoranMorgan:2006} Doran and Morgan have shown that there are 14 possible classes of variations of Hodge structure which can be associated to families of Calabi-Yau threefolds $Y$ with $h^{2,1} = 1$. They have provided explicit examples for all but one of these classes and given some suggestions on how one could construct an example for the last class (the ``14th case'' referred to in the title of this paper).

By analogy with other examples, one could hope to start with a complete intersection with $h^{1,1} = 1$ in the weighted projective space $\WP(1,1,1,1,4,6)$. Unfortunately this ambient space is not Fano, so the Batyrev-Borisov mirror construction based on nef-partitions, described in \ref{sec:NefCI}, cannot be applied to obtain a family with $h^{2,1} = 1$. Instead, Kreuzer and Sheidegger have suggested working with a slightly different ambient space, a non-crepant blow-up of $\WP(1,1,1,1,4,6)$, which is Fano and has a family of complete intersections corresponding to a suitable nef-partition, so the mirror transition is possible (See~\cite[Section~8 and Appendix E.2]{KKRS:2005} for some discussion of this example). We will construct and explore this mirror family below using Sage.

Let $\Delta \subset M_\RR$ be a $5$-dimensional reflexive polytope with polar given by%
\begin{sagecommandline}
sage: Delta5_polar = LatticePolytope([(1,-1,0,0,0), (-1,1,0,0,0), (-1,-1,0,0,0), (-1,-1,2,0,0), (12,0,-1,-1,-1), (0,12,-1,-1,-1), (0,0,-1,-1,-1), (0,0,11,-1,-1), (0,0,-1,2,-1), (0,0,-1,-1,1)], lattice=ToricLattice(5))
\end{sagecommandline}
i.e. the vertices of $\Delta^\circ$ are given by the columns of the following matrix
\begin{equation} \label{eq:Delta5 polar}
\sage{Delta5_polar.vertices_pc()}.
\end{equation}
The nef-partition we are interested in is ``in agreement'' with the block structure of this matrix: one part is formed by the first four vertices and the other by the last six.%
\begin{sagecommandline}
sage: np = NefPartition([0]*4+[1]*6, Delta5_polar)
sage: np
Nef-partition {0, 1, 2, 3} U {4, 5, 6, 7, 8, 9}
\end{sagecommandline}

Let $\Sigma$ be a crepant subdivision of $\Sigma_\Delta$ and $X = X_\Sigma$ be the corresponding crepant partial resolution of $X_\Delta$. The choice of this resolution will depend on our needs, but for the moment we only ensure that $\Sigma$ is simplicial. We will use $y_i$ to denote homogeneous coordinates on $X$, with $i$ being the index of the corresponding point of $\Delta^\circ$. For future use we also introduce parameters $c$, $d$, and $e$ into the base field of $X$. Finally, we let $Y \subset X$ be a generic member of the family of complete intersections corresponding to the nef-partition above.%
\begin{sagecommandline}
sage: X5 = CPRFanoToricVariety(np.Delta(), make_simplicial=True, coordinate_names="y+", base_field=QQ["c,d,e"].fraction_field())
sage: Y = X5.nef_complete_intersection(np)
\end{sagecommandline}
The defining polynomials of $Y$ are
\begin{align}
\label{eq:g0 general}
g_0
&=
a_{0} y_{0}^{2} y_{4}^{12} + a_{1} y_{1}^{2} y_{5}^{12} + a_{2} y_{0}y_{1}y_{2}y_{3}
,\\
\nonumber
g_1
&=
b_{4} y_{4}^{6}y_{5}^{6}y_{6}^{6}y_{7}^{6}
+ b_{5} y_{4}^{4}y_{5}^{4}y_{6}^{4}y_{7}^{4}y_{8} 
+ b_{3} y_{2}^{2}y_{6}^{12} 
+ b_{2} y_{3}^{2}y_{7}^{12} 
\\
&+ b_{7} y_{4}^{3}y_{5}^{3}y_{6}^{3}y_{7}^{3}y_{9}
+ b_{6} y_{4}^{2}y_{5}^{2}y_{6}^{2}y_{7}^{2}y_{8}^{2}
+ b_{8} y_{4}y_{5}y_{6}y_{7}y_{8}y_{9} 
+ b_{0} y_{8}^{3} 
+ b_{1} y_{9}^{2}.
\label{eq:g1 general}
\end{align}
There are 12 parameters in these equations, but their number can be significantly reduced. First of all, we can use a change of variables to set $b_5 = b_6 = b_7 = 0$. To see this more easily, we switch to an affine chart.
\begin{sagecommandline}
sage: Yap = Y.affine_patch(21)
sage: g1 = Yap.defining_polynomials()[1]
\end{sagecommandline}
In this chart $g_1$ takes the form
\begin{equation*}
g_1 = \sage{g1}
\end{equation*}
and making a substitution $y_8 = y_8 + c$, $y_9 = y_9 + d + e y_8$ does not lead to any new monomials.
\begin{sagecommandline}
sage: X5 = Y.ambient_space()
sage: X5.inject_coefficients();
sage: X5.inject_variables();
sage: g1s = g1.subs(y8=y8+c, y9=y9+d+e*y8)
sage: g1s.monomials()
[y8^3, y2^2, y3^2, y8^2, y8*y9, y9^2, y8, y9, 1]
sage: g1s.monomial_coefficient(y8^2)
e^2*b1 + 3*c*b0 + e*b8 + b6
sage: g1s.monomial_coefficient(y9)
2*d*b1 + c*b8 + b7
sage: g1s.monomial_coefficient(y8)
3*c^2*b0 + 2*d*e*b1 + c*e*b8 + 2*c*b6 + e*b7 + d*b8 + b5
\end{sagecommandline}
From this we see that one can pick $c$, $d$, and $e$ to make the coefficients of the 3 monomials $y_8$, $y_8^2$ and $y_9$ vanish, leaving only 9 parameters. Since we can also scale both polynomials and 5 of the variables, we can further reduce the number of parameters to 2. 

On the other hand, computing stringy Hodge numbers of $Y$ (one can use the generating function from~\cite{BatyrevBorisov:1996STHN} implemented in PALP~\cite{PALP} or the closed form expressions from~\cite{DoranNovoseltsev:2010}) we obtain
\begin{gather*}
h_{st}^{1,1}(Y) = 243,
\qquad
h_{st}^{2,1}(Y) = 3.
\end{gather*}
If we take $X$ to be a MPCP-desingularization of $X_\Delta$, then $Y$ is generically smooth and its stringy Hodge numbers coincide with its regular ones, so the dimension of the space of complex deformations of $Y$ is $3$. The fact that we have reduced the number of parameters in the defining polynomials of $Y$ to $2$ suggests that the dimension of the space of polynomial deformations of $Y$ is $h_{poly}^{2,1}(Y) = 2$. To confirm this observation we compute the toric part of $h_{st}^{1,1}(Y^\circ)$ for the Batyrev-Borisov mirror $Y^\circ$ of $Y$ inside a MPCP-desingularization $X^\circ$ of $X_\nabla$ (it is easier to work with a maximal resolution in this case, since $\nabla^\circ$ is much smaller than $\Delta^\circ$, i.e. it has only a few lattice points.)
\begin{sagecommandline}
sage: X5m = CPRFanoToricVariety(np.nabla(), make_simplicial=True, coordinate_points="all")
sage: X5m.is_smooth()
True
sage: Ym = X5m.nef_complete_intersection(np.dual(), monomial_points="vertices")
sage: H = Ym.ambient_space().cohomology_ring()
sage: H.gens()
([z5], [z5], [4*z5 + 2*z6 + z7 + z9], [6*z5 + 3*z6 + 2*z7 + z8 + z9], [z5], [z5], [z6], [z7], [z8], [z9])
\end{sagecommandline}
Here we have restricted the monomials used in the defining equations of $Y^\circ$ to save some time, since the nef divisors corresponding to it have hundreds of monomial sections and we only need its cohomology class. Note that $X^\circ$ is smooth, so $Y^\circ$ is (generically) smooth as well, giving $h_{st}^{1,1}(Y^\circ) = h^{1,1}(Y^\circ)$. To compute $h_{tor}^{1,1}(Y^\circ)$ (the contribution of $H^{1,1}(X^\circ)$ to $H^{1,1}(Y^\circ)$), we just need to consider the intersections of the cohomology class of $Y^\circ$ with the generators of the cohomology ring of $X^\circ$. As can be seen from the above output, this ring can be generated by the last five of ``all generators'', each of which corresponds to a torus-invariant subvariety of codimension $1$, e.g. \sagecommand{[z5]} corresponds to $\set{z_5 = 0}$, whose cohomology class is equivalent to $\set{z_0 = 0}$, $\set{z_1 = 0}$, and $\set{z_4 = 0}$.
\begin{sagecommandline}
sage: Ym_c = Ym.cohomology_class()
sage: Ym_c
[24*z5^2 - 6*z6^2 + 8/3*z7*z8 - 2/3*z8^2 - 13*z6*z9 - 3/2*z9^2]
sage: [Ym_c * g for g in H.gens()[-5:]]
[[24*z5^3 + 3*z6^3 + 14/9*z7*z8^2 + 1/9*z8^3 + 19/2*z6^2*z9 + 3/8*z9^3], [-z6^2*z9], [0], [0], [0]]
\end{sagecommandline}
From this we conclude that $h_{tor}^{1,1}(Y^\circ) = h_{poly}^{2,1}(Y) = 2$.

As pointed out in~\cite{DoranMorgan:2006}, a certain subfamily of the complete intersections $Y$ has the desired hypergeometric series corresponding to the 14th case. To describe that subfamily precisely in our setting we will compute the GKZ series of $Y$, following the algorithm outlined in \cite[Section~5.5]{CoxKatz:1999} and \cite[Appendix~A]{KKRS:2005}. 

It is easy to check that the monomials that can be eliminated in the equations for $Y$ correspond to all points of $\Delta_i$ that are neither vertices of $\Delta_i$ nor the origin (whilst this is true for our particular case, in general such information cannot be easily determined and one has to perform computations in the cohomology ring as above). Let's use this information to reconstruct our varieties:
\begin{sagecommandline}
sage: X5 = CPRFanoToricVariety(np.Delta(), make_simplicial=True, coordinate_names="y+")
sage: Y = X5.nef_complete_intersection(np, monomial_points="vertices+origin")
sage: X5m = CPRFanoToricVariety(np.nabla())
sage: X5m.Mori_cone().rays()
((0, 0, 2, 3, 0, 0, 1, -6), (1, 1, 0, 0, 1, 1, -2, -2))
\end{sagecommandline}
We now use the generators of the Mori cone of $X^\circ$ to construct moduli parameters of $Y$. These generators are given as elements of the row span of the Gale transform of the fan of $X^\circ$, the $i$-th element of each generator corresponds to the $i$-th ray of this fan, except for the last one which corresponds to the origin and is equal to the negative sum of other entries. Note that for complete intersections we need to take such sums for each part of the nef-partition separately; the following code adds them in a way that is compatible with the order of coefficients of $Y$ in Sage, with $a_2$ and $b_8$ corresponding to the origin (the coefficient-monomial correspondence for the newly constructed $Y$ is the same as in~\eqref{eq:g0 general} and~\eqref{eq:g1 general}, since the coefficient indices come from the internal enumeration of polytope lattice points).
\begin{sagecommandline}
sage: coefs = Y.ambient_space().base_ring().gens()
sage: coefs
(a0, a1, a2, b0, b1, b2, b3, b4, b8)
sage: degrees = [[[ray[i] for i in p] for p in np.dual().parts()] for ray in X5m.Mori_cone()]
sage: degrees = [flatten([p + [-sum(p)] for p in degs]) for degs in degrees]
sage: matrix(degrees)
[ 0  0  0  2  3  0  0  1 -6]
[ 1  1 -2  0  0  1  1 -2  0]
sage: B = [prod(c^d for c, d in zip(coefs, ds)) for ds in degrees]
\end{sagecommandline}

We find that the modular parameters of $Y$ are
\begin{equation*}
B_0 = \sage{B[0]},
\qquad
B_1 = \sage{B[1]},
\end{equation*}
and its GKZ series is
\begin{gather*}
\sum_{m, n}
\frac{(2m)! (6n)!}{(m!)^4 (2n)! (3n)! (n - 2m)!} B_1^m B_0^n,
\end{gather*}
where the summation is over all integers $m$ and $n$ such that the arguments of all factorials are non-negative. Making a substitution $(n - 2m) \to n$, we can sum over all non-negative integers:
\begin{gather*}
\sum_{m, n \in \ZZ_{\geqslant 0}}
\frac{(2m)! (12m + 6n)!}{(m!)^4 (4m + 2n)! (6m + 3n)! n!}
B_1^m B_0^{2m + n}.
\end{gather*}

To simplify the description of $Y$ further, we now scale its defining polynomials and coordinates to set all coefficients to $1$ except for $b_3 = \xi_0$ and $b_4 = \xi_1$.
\begin{sagecommandline}
sage: Y = X5.nef_complete_intersection(np, monomial_points= "vertices+origin",coefficients=[[1,1,1],[1,1,1,"xi0","xi1",1]])
\end{sagecommandline}
The defining polynomials of $Y$ become
\begin{align}
g_0
&=
y_{0}^{2}y_{4}^{12} + y_{1}^{2}y_{5}^{12} + y_{0}y_{1}y_{2}y_{3}
,\label{eq:g0} \\
g_1
&=
\xi_{1} y_{4}^{6}y_{5}^{6}y_{6}^{6}y_{7}^{6}
+ \xi_{0} y_{2}^{2}y_{6}^{12} 
+ y_{3}^{2}y_{7}^{12} 
+ y_{4}y_{5}y_{6}y_{7}y_{8}y_{9} 
+ y_{8}^{3} 
+ y_{9}^{2}\label{eq:g1},
\end{align}
and the GKZ series takes the form
\begin{gather*}
\sum_{m, n \in \ZZ_{\geqslant 0}}
\frac{(2m)! (12m + 6n)!}{(m!)^4 (4m + 2n)! (6m + 3n)! n!}
\xi_0^m \xi_1^n.
\end{gather*}
Comparing this series with the one given in the end of~\cite{DoranMorgan:2006}, we see that the subfamily of interest is $\xi_1 = 0$, leading to the series
\begin{gather*}
\sum_{m \in \ZZ_{\geqslant 0}}
\frac{(2m)! (12m)!}{(m!)^4 (4m)! (6m)!}
\xi_0^m.
\end{gather*}
This subfamily, henceforth denoted by $Y_1$, will be of further interest in our pursuit of a geometric model for the 14th case; we will come back to it in \ref{subsection:desingularization}. 

We complete this subsection by showing that the complete intersections $Y$ admit a natural fibration structure, which will allow us to access many of their properties. Going back to the definition of $\Delta^\circ$ in~\eqref{eq:Delta5 polar}, we observe that, in addition to its column decomposition into a nef-partition, it has a ``natural'' row decomposition. Indeed, the projection onto the first two coordinates corresponds to a toric fibration $\tld{\alpha}\colon X \to B$ over a $2$-dimensional toric variety $B$, as long as we pick a compatible resolution of $X_\Delta$. To get such a resolution, we start with the face fan of $\Delta^\circ$, take its minimal subdivision compatible with the projection, and then subdivide it to get a simplicial fan, so that $X$ is an orbifold.
\begin{sagecommandline}
sage: m = matrix([(1,0), (0,1)] + [(0,0)]*3)
sage: Delta2_polar = LatticePolytope([(1,-1), (-1,1), (-1,-1), (1,0), (0,1)])
sage: B2 = CPRFanoToricVariety(Delta_polar=Delta2_polar, coordinate_names="u+")
sage: Sigma5 = FaceFan(Delta5_polar)
sage: Sigma5.nrays(), Sigma5.ngenerating_cones()
(10, 14)
sage: Sigma5 = FanMorphism(m, Sigma5, B2.fan(), subdivide=True).domain_fan()
sage: Sigma5.nrays(), Sigma5.ngenerating_cones()
(10, 22)
sage: X5 = CPRFanoToricVariety(np.Delta(), charts=[C.ambient_ray_indices() for C in Sigma5], coordinate_names="y+", make_simplicial=True, check=False)
sage: alpha = FanMorphism(m, X5.fan(), B2.fan())
sage: alpha.is_fibration()
True
\end{sagecommandline}

In homogeneous coordinates
\begin{gather*}
\tld{\alpha}\colon
\left[y_0 : \dots : y_9\right]
\mapsto
\left[u_0 : \dots : u_4\right]
=
\left[y_{0} : y_{1} : y_{2} y_{3} : y_{4}^{12} : y_{5}^{12}\right].
\end{gather*}
Note that the hypersurface defined by the polynomial $g_0$ in~\eqref{eq:g0} depends only on the variables involved in the projection map. This means that we can interpret $g_0 = 0 $ as a defining equation of a \emph{curve} $C = \tld{\alpha}(\set{g_0 = 0}) \subset B$ and $g_1 = 0$ as a defining equation of a \emph{surface} in each fibre of $\tld{\alpha}$, in other words, $\tld{\alpha}$ induces a fibration of the complete intersection $Y$ over $C$.

Generic fibres of $\tld{\alpha}\colon X \to B$ correspond to the fan whose rays are generated by the last four vertices of $\Delta^\circ$. The polytope spanned by these vertices is the last (the 4318-th) $3$-dimensional reflexive polytope in the Kreuzer-Skarke list (included in Sage), with its normal fan corresponding to $\WP(1,1,4,6)$:
\begin{sagecommandline}
sage: vertices = [r[2:] for r in alpha.kernel_fan().rays()]
sage: p = LatticePolytope(vertices)
sage: p.is_reflexive()
True
sage: p.index()
4318
sage: NormalFan(p).rays()
(N(1, 0, 0), N(0, 1, 0), N(0, 0, 1), N(-1, -4, -6))
\end{sagecommandline}

Summarizing the results of this section, we find:

\begin{proposition} \label{prop:Yproperties} The complete intersection model $Y$ is generically a smooth Calabi-Yau threefold having $h^{2,1}(Y) = 3$ and $h^{2,1}_{poly}(Y) = 2$. There is a torically induced fibration $\tilde{\alpha}\colon Y \to C$ onto a curve $C$ whose generic fibre is an anticanonical hypersurface in the three-dimensional space polar to $\WP(1,1,4,6)$.

Furthermore, the GKZ series of the subfamily $Y_1$ obtained by setting the modular parameter $\xi_1 = 0$ coincides with the GKZ series predicted for the 14th case VHS.
\end{proposition}

\subsection{Anticanonical Hypersurfaces}
\label{subsection:AH model} 

In the previous section we were able to represent the complete intersection $Y$ as a fibration over a curve with fibres living in the space polar to $\WP(1,1,4,6)$. This suggests that the 14th case could also be realized by a family of anticanonical hypersurfaces in a four-dimensional space which can be fibred by the same toric varieties.

To try to find such a family, we searched\footnote{We performed this search using Sage, but have chosen not to include the corresponding code in the text due to its length. The algorithm used simply runs over all hyperplanes that intersect a given polytope in a sublattice polytope, with a few technical tricks to make this search happen in a reasonable timeframe and to store the results of big searches compactly. The interested reader may find a Sage worksheet implementing this search attached to the arXiv submission of this paper (note, however, that in order to run this worksheet it is necessary to install the optional package \sagecommand{polytopes_db_4d}; beware of the 8.7GB download size!)} for Fano varieties fibred by the space polar to $\WP(1,1,4,6)$ among those whose anticanonical hypersurfaces have small $h^{2,1}$, with the extra condition that the torically induced fibration is ``balanced'': this means that the same 3-dimensional reflexive polytopes can play the roles of both slices and projections, as described in \ref{subsection:torically induced fibrations}. The space polar to $\WP(1,1,2,8,12)$ satisfies these requirements; below we study the family of anticanonical hypersurfaces in it.

Let $\Delta \subset M_\RR$ be a $4$-dimensional reflexive polytope given by
\begin{sagecommandline}
sage: Delta4 = LatticePolytope([(1,0,0,0), (0,1,0,0), (0,0,1,0), (0,0,0,1), (-1,-2,-8,-12)])
\end{sagecommandline}
i.e. the vertices of $\Delta$ and of $\Delta^\circ \subset N_\RR$ are given by columns of the following matrices
\begin{equation*}
\sage{Delta4.vertices_pc()},
\qquad
\sage{Delta4.polar().vertices_pc()}.
\end{equation*}
Let $\Sigma$ be a crepant subdivision of $\Sigma_\Delta$ and $X = X_\Sigma$ be the corresponding crepant partial resolution of $X_\Delta$. As before the choice of this resolution will depend on our needs; for now we add only one extra ray in addition to the vertices of $\Delta^\circ$, specifically the ray corresponding to the midpoint $(11, -1, -1, -1)$ between the $0$-th and the $3$-rd vertices, which is the $16$-th point in the internal enumeration in Sage:
\begin{sagecommandline}
sage: Delta4.polar().point(16)
(11, -1, -1, -1)
sage: X4 = CPRFanoToricVariety(Delta=Delta4, coordinate_points=range(5)+[16])
sage: B1 = toric_varieties.P1("s,t")
sage: beta = FanMorphism(matrix(4,1,[1,1,4,6]), X4.fan(), B1.fan())
sage: beta.is_fibration()
True
\end{sagecommandline}
This extra ray is necessary to make $\Sigma$ compatible with the projection onto the line in the direction $(1,1,4,6)$. If $B = \PP^1$ with coordinates $[s:t]$ corresponding to the (unique) complete fan on this line and $\tld{\beta}\colon X \to B$ is the toric morphism associated to this projection, then
\begin{gather}
\tld{\beta}\colon
\left[z_{0} : z_{1} : z_{2} : z_{3} : z_{4} : z_{16}\right]
\mapsto
[s : t]
=
\left[z_{0}^{12} : z_{3}^{12}\right].
\label{eq:beta}
\end{gather}
Let $Z \subset X$ be a generic anticanonical hypersurface. Its defining polynomial is
\begin{multline*}
h =
a_{0} z_{0}^{24}z_{16}^{12} 
+ a_{5} z_{0}^{12}z_{3}^{12}z_{16}^{12} 
+ a_{4} z_{3}^{24}z_{16}^{12} 
+ a_{6} z_{0}^{6}z_{2}^{6}z_{3}^{6}z_{16}^{6} \\
+ a_{1} z_{2}^{12}
+ a_{10} z_{0}z_{1}z_{2}z_{3}z_{4}z_{16} 
+ a_{2} z_{1}^{3} 
+ a_{3} z_{4}^{2}.
\end{multline*}
Scaling the whole polynomial and four independent coordinates we can eliminate 5 out of 8 parameters. Using Batyrev's formulas for the Hodge numbers of anticanonical hypersurfaces we check that $h^{2,1}(Z) = h^{2,1}_{poly}(Z) = 3$, so we should indeed have 3 independent parameters. Summarizing, we find:

\begin{proposition} \label{prop:Zproperties} The hypersurface model $Z$ is generically a smooth Calabi-Yau threefold with $h^{2,1}(Z) = h^{2,1}_{poly}(Z)  = 3$. There is a torically induced fibration $\tilde{\beta}\colon Z \to B \cong \PP^1$ whose general fibre is an anticanonical hypersurface in the three-dimensional space polar to $\WP(1,1,4,6)$.
\end{proposition}

Anticanonical hypersurfaces inside the space polar to $\WP(1,1,2,8,12)$ were extensively studied by Bill\'{o} et al.~\cite{Billo:1998}; in order to conveniently use their results we will match our toric description with theirs (they considered hypersurfaces in $\WP(1,1,2,8,12)$ with extra symmetries, which allow taking the quotient under a certain group action).

First, we rewrite the polynomial of $Z$ in its ``fibred'' form, thinking of it as a polynomial in $z_1, z_2, z_4, z_{16}$ only and working in a chart with $t = z_3 = 1$:
\begin{align*}
\left[a_{0} s^2 + a_{5} s + a_{4}\right] z_{16}^{12} 
+ \left[a_{6} z_{0}^{6}\right] z_{2}^{6}  z_{16}^{6} 
+ a_{1} z_{2}^{12} 
+ \left[a_{10} z_{0}\right] z_{1} z_{2} z_{4} z_{16} 
+ a_{2} z_{1}^{3} 
+ a_{3} z_{4}^{2}.
\end{align*}
Now we can compare our representation with equation (4.19) in~\cite{Billo:1998}:
\begin{gather*}
W^{(2)}(x; B', \psi_0, \psi_1) 
=
\frac{1}{12} (B' x_0^{12} + x_3^{12}) + \frac{1}{3} x_4^3 + \frac{1}{2} x_5^2 - \psi_0 x_0 x_3 x_4 x_5 - \frac{1}{6} \psi_1 x_0^6 x_3^6.
\end{gather*}
We see that the matching of coordinates and coefficients is
\begin{align*}
z_1 &= x_4, &
z_2 &= x_3, &
z_4 &= x_5, &
z_0 z_{16} &= x_0,\\
a_1 &= \frac{1}{12}, &
a_2 &= \frac{1}{3}, &
a_3 &= \frac{1}{2}, &
a_6 &= -\frac{1}{6} \psi_1, &
a_{10} &= -\psi_0, 
\end{align*}
and
\begin{gather*}
a_0 s + \frac{a_4}{s} + a_5 = \frac{1}{12} B'.
\end{gather*}
To match these remaining parameters we use the definition of $B'$ given by equation (3.18) in~\cite{Billo:1998}:
\begin{gather*}
B' = \frac{1}{2}\left(B\zeta + \frac{B}{\zeta} - 2\psi_s\right),
\end{gather*}
where $\zeta$ is an affine coordinate on the base of the fibration, so
\begin{align*}
s &= \zeta, &
a_0 &= \frac{B}{24}, &
a_4 &= \frac{B}{24}, &
a_5 &= -\frac{\psi_s}{12}.
\end{align*}
We can use these parameters in Sage as follows:
\begin{sagecommandline}
sage: var("B,psi0,psi1,psi_s");
sage: Z = X4.anticanonical_hypersurface( coefficients=[B/24,1/12,1/3,1/2,B/24,-psi_s/12,-psi1/6,-psi0])
\end{sagecommandline}
Now the defining polynomial of $Z$ has the form
\begin{multline*}
\frac{B}{24} z_{0}^{24}z_{16}^{12} 
- \frac{\psi_s}{12} z_{0}^{12}z_{3}^{12}z_{16}^{12} 
+ \frac{B}{24} z_{3}^{24}z_{16}^{12} 
- \frac{\psi_1}{6} z_{0}^{6}z_{2}^{6}z_{3}^{6}z_{16}^{6} 
\\
+ \frac{1}{12} z_{2}^{12} 
- \psi_{0} z_{0}z_{1}z_{2}z_{3}z_{4}z_{16} 
+ \frac{1}{3} z_{1}^{3} 
+ \frac{1}{2} z_{4}^{2}.
\end{multline*}
One of the four parameters in this representation is redundant, e.g. we can set $\psi_0 = 1$.

\subsection{Geometric Transitions}

Both the complete intersections $Y$ and anticanonical hypersurfaces $Z$ presented above have $h^{2,1} = 3$, while the original goal in their construction was to obtain families with $h^{2,1} = 1$ that could provide geometric examples of the 14th case. In fact, in the hypersurface case it was known in advance that we would ``fail'', since it is known that there are only 5 reflexive polytopes yielding Calabi-Yau threefolds with $h^{2,1} = 1$ and \cite{DoranMorgan:2006} showed that they already provide examples for other classes of Hodge structure variations.

However, we can still try to obtain our desired families by using subfamilies of the constructed ones. Of course, simply fixing some of the parameters does not change the Hodge numbers of the Calabi-Yau threefolds in question, but if these threefolds were to become singular we could try to resolve the singularities and hope that $h^{2,1} = 1$ holds for the resolved family. So we are now looking for geometric transitions from the already constructed families to some new ones with, hopefully, ``correct'' Hodge numbers. 

The existence of such geometric transitions was first suggested by Clemens \cite{Clemens:1983,Clemens:1983a} and later expanded upon by Friedman \cite{Friedman:1986}. In this subsection we give a precise definition and a basic classification of geometric transitions, following~\cite{Rossi:2010}.

\begin{definition}
Let $Y$ and $\tld{Y}$ be smooth Calabi-Yau threefolds. They are connected by a \keyterm{geometric transition} if there exist a normal variety $\overline{Y}$, a birational contraction $\phi\colon Y \to \overline{Y}$, and a complex deformation (smoothing) of $\overline{Y}$ to $\tld{Y}$.
It is a \keyterm{primitive geometric transition} if $\phi$ cannot be factored into birational morphisms of normal varieties.
It is a \keyterm{conifold transition} if $\overline{Y}$ has only conifold singularities (ordinary double points).
It is a \keyterm{trivial geometric transition} if $\tld{Y}$ is a deformation of $Y$.
\end{definition}

\begin{theorem}
Let $\phi\colon Y \to \overline{Y}$ be a primitive contraction of a smooth Calabi-Yau threefold $Y$ to a normal variety $\overline{Y}$. Let $E$ be the exceptional locus of $\phi$. Then $\phi$ is of one of the following three types:
\begin{enumerate}
\item[\keyterm{Type I}] $\phi$ is small, $E$ may be reducible and is composed of finitely many rational curves;
\item[\keyterm{Type II}] $\phi$ contracts a divisor to a point, $E$ is irreducible and is a generalized del Pezzo surface;
\item[\keyterm{Type III}] $\phi$ contracts a divisor to a smooth curve $C$, $E$ is irreducible and is a conic bundle over $C$.
\end{enumerate}
\end{theorem}

\begin{proof}
See~\cite{Wilson:1992, Wilson:1993}, the given formulation is \cite[Theorem~1.9]{Rossi:2010}.
\end{proof}

\begin{definition}
\label{def:primitive transitions}
A primitive geometric transition is of type I, II, or III, if the corresponding birational contraction is of type I, II, or III, respectively.
\end{definition}

The existence of geometric transitions relating Calabi-Yau threefolds with different Hodge numbers has interesting consequences for their moduli spaces. In particular, it has been suggested \cite{Reid:1987, Friedman:1991} that geometric transitions could be used to connect together the moduli spaces of all compact complex threefolds with trivial canonical bundle. Moreover, it is also natural to ask how geometric transitions interact with the string-theoretic phenomenon of mirror symmetry \cite{Dixon:1988,CandelasLynkerSchimmrigk:1990,GreenePlesser:1990} and this has been the object of a great deal of study; a good summary of this theory is given by Morrison \cite{Morrison:1999}. We will return to the question of mirror symmetry in our context in  \ref{section:mirrormoduli}.

However, before we use geometric transitions to construct new families, we first find that our two models $Y$ and $Z$ are in fact connected via a type III geometric transition. To see this explicitly, we begin by comparing their fibration structures.

\section{Matching the Models} \label{sect:matching}
\subsection{K3 Fibrations}
\label{subsection:K3 fibrations}

The generic fibres of the fibrations induced by $\tld{\alpha}$ and $\tld{\beta}$ constructed earlier are anticanonical hypersurfaces inside the three-dimensional space polar to $\WP(1,1,4,6)$, i.e. they are generically two-dimensional Calabi-Yau varieties: K3 surfaces. The choice of toric ambient space for them induces a lattice polarization, in this case by the lattice $M = H \oplus E_8 \oplus E_8$, where $H$ is the hyperbolic lattice of rank $2$ and $E_8$ is the unique even negative-definite unimodular lattice of rank $8$. Such $M$-polarized K3 surfaces (the name ``$M$-polarized'' is a bit unfortunate in the toric context, but it should not cause too much confusion) were originally studied by Shioda and Inose \cite{ShiodaInose:1977, Inose:1978} and have recently been revisited in \cite{shioda:2006}, \cite{ClingherDoran:2007} and~\cite{CDLW:2009}, we start this section with a summary of their properties.

\begin{definition} \label{definition:M-polarized}
An \keyterm{$M$-polarization} on a K3 surface $X$ is a primitive lattice embedding $i\colon M \hookrightarrow \NS(X)$, such that the image $i(M)$ in the N\'eron-Severi lattice $\NS(X)$ contains a pseudo-ample class (corresponding to an effective nef divisor with positive self-intersection).
\end{definition}

\begin{definition} 
Let $X$ be a K3 surface. An involution $i$ on $X$ is a \keyterm{Nikulin involution} if $i^* \omega = \omega$ for any holomorphic 2-form $\omega$ on $X$.
\end{definition}

\begin{theorem}
\label{theorem:M-polarized K3}
Let $X$ be an $M$-polarized K3 surface. Then
\begin{enumerate}
\item
$X$ is isomorphic to the minimal resolution of a quartic surface in $\PP^3$ given by
\begin{gather*}
y^2 z w - 4 x^3 z + 3 a x z w^2 + b z w^3 - \frac{1}{2} (d z^2 w^2 + w^4) = 0;
\end{gather*}
\item
the parameters $a$, $b$, and $d$ in the above equation specify a unique point $(a,b,d) \in \WP(2,3,6)$ with $d \neq 0$;
\item
$X$ canonically corresponds to a pair of elliptic curves $\set{E_1, E_2}$;
\item
the modular parameters of $X$ and $\set{E_1, E_2}$ are related by
\begin{gather*}
\pi = j(E_1) j(E_2) = \frac{a^3}{d}
\qquad \text{and} \qquad
\sigma = j(E_1) + j(E_2) = \frac{a^3 - b^2 + d}{d};
\end{gather*}
\item
generically there are exactly two isomorphism classes of elliptic fibrations with section on $X$: the ``standard'' fibration, which has an $H \oplus E_8 \oplus E_8$ polarization realized via two type $II^*$ exceptional fibres, and the ``alternate'' fibration which has an $H \oplus D_{16}$ polarization realized via a type $I_{12}^*$ exceptional fibre;
\item
there exists a Nikulin involution on $X$.
\end{enumerate}
\end{theorem}

\begin{proof}
See \cite{shioda:2006}, \cite[Theorem~1.1, Corollary~1.3, Section~3]{ClingherDoran:2007} and \cite[Theorems~3.1,~3.2]{CDLW:2009}.
\end{proof}

\begin{proposition}
\label{prop:toric M-polarized K3}
An anticanonical hypersurface in the space polar to the weighted projective $\WP(1,1,4,6)$ is an $M$-polarized K3 surface defined by
\begin{sagesilent}
Delta3 = LatticePolytope([(1,0,0), (0,1,0), (0,0,1), (-1,-4,-6)])
X3 = CPRFanoToricVariety(Delta3, coordinate_names="x0 x3 x1 x2")
K3 = X3.anticanonical_hypersurface(coefficient_names="lambda0, lambda3, lambda2, lambda1, lambda4, lambda5")
\end{sagesilent}
\begin{equation*}
\sage{K3.defining_polynomials()[0]} = 0.
\end{equation*}
It is related to the normal form given in Theorem~\ref{theorem:M-polarized K3} by
\begin{gather*}
a^3 = \frac{1}{12^6 \Lambda_0^2 \Lambda_1},\ 
b^2 = \frac{(6 \cdot 12^2 \Lambda_0 - 1)^2}{12^6 \Lambda_0^2 \Lambda_1},\
d = 1,
\text{ with }
\Lambda_0 = \frac{\lambda_2^3 \lambda_3^2 \lambda_4}{\lambda_5^6},\ 
\Lambda_1 = \frac{\lambda_0 \lambda_1}{\lambda_4^2}.
\end{gather*}
\end{proposition}

\begin{proof}
See \cite[Section~3.4]{CDLW:2009}\footnote{There was a typo in the preprint version of~\cite{CDLW:2009} posted on arXiv: the numerator of the expression for $b^2$ was not squared.}\!.
\end{proof}

Both the standard and alternate elliptic fibrations have toric realizations on the anticanonical hypersurfaces in the above proposition. Here we will point out a few key facts about this toric picture that will be used later; for a detailed treatment of the combinatorics of torically induced elliptic fibrations of K3 surfaces we refer the reader to~\cite{PerevalovSkarke:1997} or~\cite{Rohsiepe:2004}. 

\begin{figure}
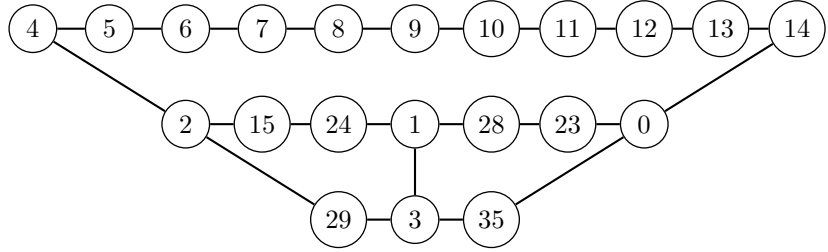

\centering
\sage{G}
\caption{Intersection of toric divisors on an \mbox{$M$-polarized} K3 surface.} \label{figure:G}
\end{figure}

\begin{figure}
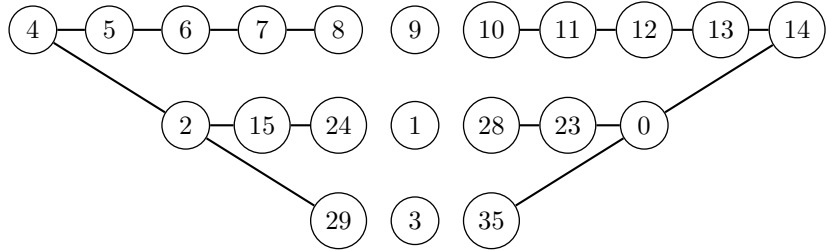

\centering
\sage{E8E8}
\caption{ADE diagram of the standard fibration on an \mbox{$M$-polarized} K3 surface.} \label{figure:E8E8}
\end{figure}

\begin{figure}
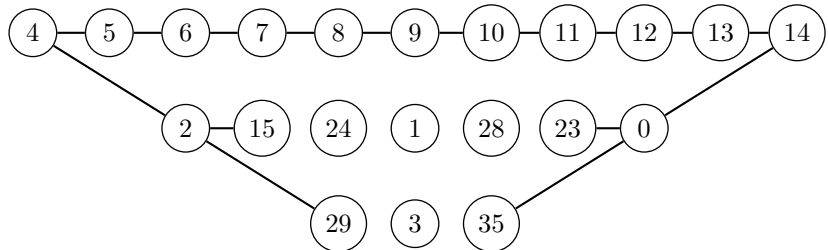

\centering
\sage{D16}
\caption{ADE diagram of the alternate fibration on an \mbox{$M$-polarized} K3 surface.} \label{figure:D16}
\end{figure}

In our case the toric correction term is zero, meaning that the intersection pattern of divisors on a generic anticanonical hypersurface, which have three linear relations between them, is given by the $1$-skeleton of the reflexive polytope in the lattice $N$. Two directions in $N$ determine the two fibrations over $\PP^1$, and the ADE diagrams of the fibres over zero and infinity are given by the parts of the $1$-skeleton that are projected into the interiors of rays of the $\PP^1$ fan. Here is an explicit construction of these diagrams:
\begin{sagecommandline}
sage: Delta3 = LatticePolytope([(1,0,0), (0,1,0), (0,0,1), (-1,-4,-6)])
sage: P = Delta3.polar()
sage: sum(e.ninterior_points() * ed.ninterior_points() for e, ed in zip(Delta3.edges(), P.edges()))
sage: G = P.skeleton()
sage: d_std = vector([1,2,3])
sage: E8E8 = G.subgraph(edge_property=lambda e: (P.point(e[0])*d_std) * (P.point(e[1])*d_std) > 0)
sage: d_alt = vector([0,1,1])
sage: D16 = G.subgraph(edge_property=lambda e: (P.point(e[0])*d_alt) * (P.point(e[1])*d_alt) > 0)
\end{sagecommandline}
\begin{sagesilent}
pos={
4: (-5, 1),
5: (-4, 1),
6: (-3, 1),
7: (-2, 1),
8: (-1, 1),
9: (0, 1),
10: (1, 1),
11: (2, 1),
12: (3, 1),
13: (4, 1),
14: (5, 1),

2: (-3, 0),
15: (-2, 0),
24: (-1, 0),
1: (0, 0),
28: (1, 0),
23: (2, 0),
0: (3, 0),

29: (-1, -1),
3: (0, -1),
35: (1, -1),
}
G.set_pos(pos)
E8E8.set_pos(pos)
D16.set_pos(pos)
G.set_latex_options(tkz_style="Normal", units="in", graphic_size=(4,1))
E8E8.set_latex_options(tkz_style="Normal", units="in", graphic_size=(4,1))
D16.set_latex_options(tkz_style="Normal", units="in", graphic_size=(4,1))
\end{sagesilent}
The graphs constructed above are shown in \vrefrange{figure:G}{figure:D16}. In these graphs, the labels correspond to the indices of lattice points of the polytope~\sagecommand{P}. For the standard fibration the toric divisor corresponding to point $9$ gives a uniquely defined section, and for the alternate fibration there are two toric sections corresponding to points $24$ and $28$. 


Next we will see how this theory can be applied to our models $Y$ and $Z$. For the complete intersection model $Y$ described in \ref{subsection:CI model} the role of $[x_0 : x_1 : x_2 : x_3]$ in Proposition~\ref{prop:toric M-polarized K3} is played by $[y_6 : y_7 : y_9 : y_8]$. Using this we obtain the following expressions for the modular parameters of a fibre:
\begin{gather*}
\pi_Y = \frac{y_{4}^{12} y_{5}^{12}}{12^6 \xi_{0} y_{2}^{2} y_{3}^{2}},
\qquad
\sigma_Y = 1 + \left(\xi_1 - 3 \cdot 12^2 \xi_{1}^2\right) \frac{y_{4}^{12} y_{5}^{12}}{12^3 \xi_{0} y_{2}^{2} y_{3}^{2}}.
\end{gather*}

\begin{figure}[ht]
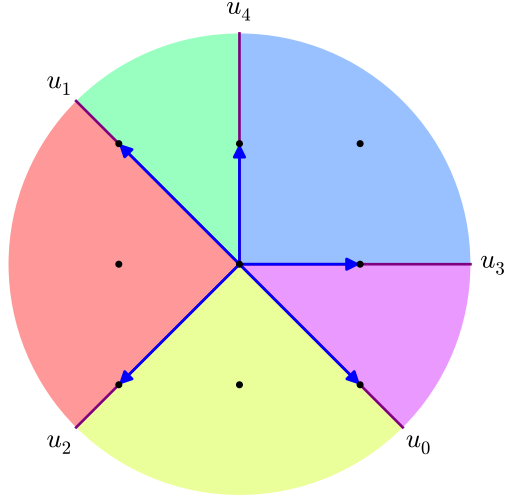

\centering
\sageplot[]{B2.plot(wall_label=None), figsize=[3, 3]}
\caption{Base $B^2$ of the fibration $\tld{\alpha}$} \label{figure:B2}
\end{figure}

For the anticanonical hypersurface model $Z$ from \ref{subsection:AH model}, $[x_0 : x_1 : x_2 : x_3]$ correspond to $[z_{16} : z_2 : z_4 : z_1]$ and we compute
\begin{align*}
\pi_Z &= \frac{\psi_0^{12}}{2} \cdot
\frac{z_0^{12} z_3^{12}}{B z_0^{24} - 2 \psi_{s} z_0^{12} z_3^{12} + B z_3^{24}},\\
\sigma_Z &= 1 - 2 (\psi_0^6 \psi_1 + \psi_1^2)
\frac{z_0^{12} z_3^{12}}{B z_0^{24} - 2 \psi_{s} z_0^{12} z_3^{12} + B z_3^{24}}.
\end{align*}
Next we switch to coordinates on the bases $B^2$ and $B^1$ of the fibrations $\tld{\alpha}$ and $\tld{\beta}$ respectively:\footnote{We have used the same notation $\Delta$, $\Sigma$, $X$, and $B$ in \ref{subsection:CI model} and \ref{subsection:AH model} to refer to different objects. When we need to consider them together and it is necessary to distinguish them, we will use their dimensions as superscripts. Names in Sage examples always include these dimensions to allow reusing of objects in later sections.}
\begin{align}
\pi_Y &= \frac{u_3 u_4}{12^6 \xi_0 u_2^2},\nonumber\displaybreak[0]\\
\sigma_Y &= 1 + \left(\xi_1 - 3 \cdot 12^2 \xi_1^2\right) \frac{u_3 u_4}{12^3 \xi_0 u_2^2},\nonumber\displaybreak[0]\\
\pi_Z &= \frac{\psi_0^{12}}{2} \cdot
\frac{s t}{B s^2 - 2 \psi_{s} s t + B t^2},\label{eq:pi_Z}\\
\sigma_Z &= 1 - 2 (\psi_0^6 \psi_1 + \psi_1^2)
\frac{s t}{B s^2 - 2 \psi_{s} s t + B t^2}.\label{eq:sigma_Z}
\end{align}
Finally, recall that $Y$ is fibred not over $B^2$, but over a curve $C \subset B^2$ corresponding to $g_0$. It is easy to see from~\eqref{eq:g0}, that $g_0$ is the pullback of
\begin{gather}
u_0^2 u_3 + u_1^2 u_4 + u_0 u_1 u_2. \label{eq:g0 on B}
\end{gather}
Using the fan of $B^2$, shown in \vref{figure:B2}, we can see that $C$ does not intersect the divisors corresponding to $u_0$ and $u_1$, e.g. if $u_0 = 0$ then $u_1^2 u_4 \neq 0$, since there are no cones containing rays corresponding to $u_0$ and $u_1$ or $u_4$.

This means that $C$ is isomorphic to the curve in $\PP^2$ with coordinates $[u : v : w] = [u_3 : u_4 : u_2]$, given by $u + v + w = 0$, i.e. $C \simeq \PP^1$ with coordinates $[u : v]$.
In terms of these new coordinates the modular parameters of a fibre are
\begin{align}
\pi_Y &= \frac{1}{12^6 \xi_0} \cdot \frac{u v}{(u + v)^2},\label{eq:pi_Y}\\
\sigma_Y &= 1 + \frac{\xi_1 - 3 \cdot 12^2 \xi_1^2}{12^3 \xi_0} \cdot  \frac{u v}{(u + v)^2}.\label{eq:sigma_Y}
\end{align}

Comparing expressions for $\pi$ and $\sigma$ for the complete intersection model $Y$, \eqref{eq:pi_Y} and \eqref{eq:sigma_Y}, with analogous expressions for the anticanonical hypersurface model $Z$, \eqref{eq:pi_Z} and \eqref{eq:sigma_Z}, we see that they are quite similar. This similarity is even more striking if we let $\psi_s = - B$ in \eqref{eq:pi_Z} and \eqref{eq:sigma_Z}, to obtain:
\begin{align}
\pi_Z &= \frac{\psi_0^{12}}{2 B} \cdot
\frac{s t}{(s + t)^2},\label{eq:pi_Z_2}\\
\sigma_Z &= 1 - \frac{2 (\psi_0^6 \psi_1 + \psi_1^2)}{B} \cdot
\frac{s t}{(s + t)^2}.\label{eq:sigma_Z_2}
\end{align}
In fact, the expressions \eqref{eq:pi_Y} and \eqref{eq:sigma_Y} and the expressions \eqref{eq:pi_Z_2} and \eqref{eq:sigma_Z_2} would be exactly the same if we let $[u : v] = [s : t]$ and
\begin{gather}
\xi_0 = \frac{2 B}{\left(12 \psi_0^2\right)^6},
\qquad
\xi_1 = -\frac{4 \psi_1}{\left(12 \psi_0^2\right)^3}.
\label{eq:Y-Z parameter matching}
\end{gather}

Finally, we note also that with respect to the fibre modular parameters $\pi_Y$ and $\sigma_Y$, the subfamily of complete intersections with $\xi_1 = 0$ is very special: in this case $\sigma_Y \equiv 1$, so all K3 fibres correspond to ``complementary'' elliptic curves with $j(E_1) + j(E_2) = 1$.

\subsection{Ambient Space Morphism}
\label{subsection:ambient space morphism}

Given such a perfect matching of modular parameters, we may try to construct a morphism between the original varieties $Y \to Z$, perhaps by constructing a morphism (ideally, a fibration) between their ambient spaces $X^5 \to X^4$. Combining all correspondences used so far, we get
\begin{gather*}
[y_0 : \dots : y_9]
\leftrightsquigarrow
\left([u : v], [x_0 : \dots : x_3]\right)
\leftrightsquigarrow
[z_0 : \dots : z_{16}],
\end{gather*}
or, more concretely,
\begin{gather*}
[y_4 : y_8 : y_7 : y_5 : y_9 : y_6]
\leftrightsquigarrow
[z_0 : z_1 : z_2 : z_3 : z_4 : z_{16}].
\end{gather*}
This correspondence does not give us a morphism between toric varieties, since only some of the homogeneous coordinates on $X^5$ are used and they are unlikely to be well-defined due to relations involving the other coordinates. But we may look for a fan morphism that will at least map the rays corresponding to these selected $y_i$ to the rays corresponding to the indicated $z_j$.
\begin{sagecommandline}
sage: m_y = Delta5_polar.vertices_pc()(4,8,7,5,9,6).matrix()
sage: m_z = Delta4.polar().points_pc()(0,1,2,3,4,16).matrix()
sage: m = m_y.solve_right(m_z)
\end{sagecommandline}
This computation shows that if such a fan morphism exists, it must be given by the matrix~\sagecommand{m}
\begin{equation*}
\sage{m},
\end{equation*}
which is unique since \sagecommand{m_y} has maximal rank. If this matrix defines a fibration, then rays of $\Sigma^5$ must be mapped onto rays of $\Sigma^4$ or the origin.
\begin{sagecommandline}
sage: rays = [r*m for r in Sigma5.rays() if not (r*m).is_zero()]
sage: from sage.geometry.cone import normalize_rays
sage: rays = set(normalize_rays(rays, None))
sage: points = sorted(Delta4.polar().points_pc().index(r) for r in rays)
sage: points
[0, 1, 2, 3, 4, 16, 331, 333, 334]
\end{sagecommandline}
Since all necessary image rays can be generated by lattice points of $\left(\Delta^4\right)^\circ$, we can pick $\Sigma^4$ to be a crepant subdivision of $\Sigma_{\Delta^4}$. One can then subdivide the current $\Sigma^5$ (the subdivision of $\Sigma_{\Delta^5}$ compatible with the fibration $\alpha$) to obtain a fibration over $X^4$, but, unfortunately, this subdivision of $\Sigma_{\Delta^5}$ will not be crepant.

To fix this problem, we go back to the defining polynomial~\eqref{eq:g0 on B} of the curve $C \subset B^2$ and the fan of $B^2$, shown in \vref{figure:B2}. We have already established that $C$ does not intersect the divisors corresponding to $u_0$ and $u_1$. Now we note that it also does not contain the point with $u_3 = u_4 = 0$, since this would imply that one of the other coordinates is zero, which is not possible. This means that $C$ is completely contained within the part of $B^2$ corresponding to the fan generated by the rays of $u_2$, $u_3$, and $u_4$. But then $Y \subset X^5$ is completely contained within the toric variety corresponding to $\Sigma^5$ \emph{without} all cones that contain rays corresponding to either $y_0$, or $y_1$, or $y_4$ and $y_5$ together. Starting with such a subfan $\Sigma^5_\mpart$ of $\Sigma^5$, we can find a subdivision that leads to a fibration $\tld{\Phi}\colon X^5_\mpart \to X^4$ with affine lines as generic fibres.
\begin{sagecommandline}
sage: selected = []
sage: for sigma in flatten(Sigma5.cones()):
...       indices = sigma.ambient_ray_indices()
...       if (0 in indices or
...           1 in indices or
...           4 in indices and 5 in indices):
...           continue
...       selected.append(sigma)
sage: Sigma5_part = Fan(cones=selected, rays=Sigma5.rays(), discard_faces=True)
sage: rays = [ray*m for ray in Sigma5_part.rays() if not (ray*m).is_zero()]
sage: rays = set(normalize_rays(rays, None))
sage: points4 = sorted(Delta4.polar().points_pc().index(r) for r in rays)
sage: points4
[0, 1, 2, 3, 4, 16, 334]
sage: Delta4.polar().point(334)
(-1, 1, 0, 0)
sage: X4 = CPRFanoToricVariety(Delta4, coordinate_points=points4)
sage: Phi = FanMorphism(m, Sigma5_part, X4.fan(), subdivide=True)
sage: Phi.is_fibration()
True
sage: Phi.kernel_fan().rays()
(N(-1, -1, 0, 0, 0),)
sage: all(ray in Delta5_polar.points_pc() for ray in Phi.domain_fan().rays())
True
sage: points5 = [Delta5_polar.points_pc().index(ray) for ray in Phi.domain_fan().rays()]
sage: points5
[2, 3, 4, 5, 6, 7, 8, 9, 752]
sage: Delta5_polar.point(752)
(0, 0, 1, 0, 0)
sage: Phi.domain_fan().is_simplicial()
True
sage: X5_part = ToricVariety(Phi.domain_fan(), coordinate_names="y+", coordinate_indices=points5)
\end{sagecommandline}
Now not only can all necessary rays in the codomain be generated by lattice points of $\left(\Delta^4\right)^\circ$, but also all necessary rays in the domain can be generated by lattice points of $\left(\Delta^5\right)^\circ$. We can therefore think of $X^5_\mpart$ as an open subset in some crepant partial desingularization of $X_{\Delta^5}$, whilst $Y \subset X^5_\mpart$ is still a Calabi-Yau variety.

From the computational point of view, we can no longer use the framework of CPR-Fano toric varieties in Sage to represent $X^5_\mpart$; instead we have to use generic toric varieties. We will however still use the same convention for naming coordinates: the missing $y_0$ and $y_1$ now reflect the fact that we are interested only in charts where they are non-zero, in which case both can be set to $1$ using relations between coordinates.

In homogeneous coordinates we get
\begin{align*}
\tld{\Phi}\colon
[y_2 : \dots : y_9 : y_{752}]
&\mapsto
[z_0 : z_1 : z_2 : z_3 : z_4 : z_{16} : z_{334}]\\
&=
[y_4 : y_8 : y_7 : y_5 : y_9 : y_6 : y_3^2 y_{752}],
\end{align*}
which also does not involve $y_2$, since it corresponds to the only ray of the kernel fan. This map can be used to pull back $Z$ and compare it with $Y$. We delay this comparison, however, since for analysis of singularities it is convenient to perform a few more subdivisions of the underlying fans first.

\section{Singular Subfamilies} \label{sect:singular}
\subsection{Hypersurfaces}
\label{subsection:singular hypersurfaces}

Generically, complete intersections $Y$ and anticanonical hypersurfaces $Z$ in MPCP-desingularizations of their ambient spaces are smooth, since the singular locus in this case has codimension at least four. However, in this paper we are also interested in the following subfamilies:
\begin{compactenum}
\item
complete intersections $Y$ with $\xi_1 = 0$, corresponding to the subfamily with the desired GKZ series, we will denote a generic member of this subfamily as $Y_1$ to emphasize dependence on a single parameter only;
\item
hypersurfaces $Z$ with $\psi_s = -B$, corresponding to the subfamily whose K3 fibration can be ``perfectly matched'' with the K3 fibration of complete intersections $Y$, we will denote a generic member of this subfamily as $Z_2$;
\item
hypersurfaces $Z$ with $\psi_s = -B$ and $\psi_1 = 0$, corresponding to the subfamily of the above subfamily with an analogue of the $\xi_1 = 0$ restriction, we will denote a generic member of this subfamily as $Z_1$.
\end{compactenum}
Generic members of these subfamilies may be singular, we analyze their singularities using results of~\cite{Billo:1998} and computer software (Sage interfacing with Magma \cite{Magma} for computing the singular loci of affine varieties). We begin in this subsection with an analysis of the singularities of the hypersurfaces $Z_2$ and $Z_1$.

The singular locus in the moduli space of hypersurfaces $Z$ is summarized in relations~(4.39) of~\cite{Billo:1998}:
\begin{align*}
S_{a1}^\pm &\colon & (\psi_0^6 + \psi_1)^2 + \psi_s &= \pm B,\\
S_{a2}^\pm &\colon & \psi_1^2 + \psi_s &= \pm B,\\
S_b^\pm &\colon & \psi_s &= \pm B,\\
S_0 &\colon & 0 &= \phantom{\pm} B.
\end{align*}
Here relations $S_{a1}^\pm$ and $S_{a2}^\pm$ are actually the same, in the sense that they are switched by an appropriate change of coordinates (there is a finite group action on the simplified polynomial moduli space of hypersurfaces). We see that the hypersurfaces $Z_2$ are singular (and, therefore, cannot be isomorphic to the full two-parameter family of complete intersections $Y$), since the condition $S_b^-$ is satisfied. The condition $\psi_1 = 0$ makes relations $S_{a2}^\pm$ and $S_b^\pm$ the same, but does not impose singularities on its own.

An explicit chart-by-chart check for singularities of $Z_2$ using our current fan $\Sigma^4$ and ignoring the orbifold structure of $X^4$ (i.e. ignoring singularities of $Z_2$ which are inevitable due to the ambient space structure) reveals a singular locus of dimension 1 in all charts involving $z_{334}$. To study this singular locus, it is convenient to put it into smooth charts of $X^4$; in order to do this we need to subdivide $\Sigma^4$ further.

In fact, when doing this one finds that it is better to start afresh with $\Sigma_{\Delta^4}$. Recall that adding the ray corresponding to the 16-th point $v_{16} = \sage{Delta4.polar().point(16)}$ of $\left(\Delta^4\right)^\circ$, which is the midpoint between the vertices $v_0 = \sage{Delta4.polar().point(0)}$ and $v_3 = \sage{Delta4.polar().point(3)}$, was necessary to ensure compatibility with the fibration $\tld{\beta}\colon X^4 \to B^1 \simeq \PP^1$. We also had to add $v_{334} = \sage{Delta4.polar().point(334)}$ to allow for a fibration $\tld{\Phi}\colon X^5_\mpart \to X^4$. This is the only interior point of the triangular face on $v_1 = \sage{Delta4.polar().point(1)}$, $v_2 = \sage{Delta4.polar().point(2)}$, and $v_4 = \sage{Delta4.polar().point(4)}$. The face fan of this triangle (in the spanned affine sublattice with $v_{334}$ being the origin) is the fan of $\WP(1,2,3)$, and the three interior points of its edges correspond to $v_{251} = \sage{Delta4.polar().point(251)}$, $v_{276} = \sage{Delta4.polar().point(276)}$, and $v_{325} = \sage{Delta4.polar().point(325)}$. Adding all these rays is sufficient to resolve all singularities of the 3-dimensional cone on $v_1$, $v_2$, and $v_4$, but the resulting subcones are still faces of singular 4-dimensional cones. This is also reflected in the homogeneous coordinate representation of $\tld{\beta}$~\eqref{eq:beta}: the coordinates on the base $B^1$ correspond to the \mbox{12th} powers of coordinates on $X^4$ (yet the defining equation of $Z$ involves first powers of all variables). The problem is that $v_0$ and $v_3$ are ``too far away'' from the slice hyperplane defining the projection to $B^1$. We can remedy the situation by adding two more points right ``above'' and ``below'' the face on $v_1$, $v_2$, and $v_4$ (which is completely contained in the slice hyperplane), namely $v_{168} = \sage{Delta4.polar().point(168)}$ and $v_{170} = \sage{Delta4.polar().point(170)}$.

Using consecutive star-like subdivisions (which are used in Sage for automatic insertion of rays), it turns out that the best sequence is the following:
\begin{compactenum}
\item add $v_{16}$ to allow the fibration $\tld{\beta}\colon X^4 \to B^1$;
\item add $v_{168}$ and $v_{170}$ to ``improve'' this fibration;
\item add $v_{334}$ to allow the fibration $\tld{\Phi}\colon X^5_\mpart \to X^4$;
\item add $v_{251}$, $v_{276}$, and $v_{325}$ to cover the divisor of $z_{334}$ by smooth charts.
\end{compactenum}
\begin{sagecommandline}
sage: X4 = CPRFanoToricVariety(Delta4, coordinate_points=[0,1,2,3,4,16,168,170,334,251,276,325])
sage: beta = FanMorphism(matrix(4,1,[1,1,4,6]), X4.fan(), B1.fan())
sage: beta.is_fibration()
True
sage: all(sigma.is_smooth() for sigma in X4.fan() if (-1,1,0,0) in sigma)
True
sage: Z = X4.anticanonical_hypersurface( coefficients=[B/24,1/12,1/3,1/2,B/24,-psi_s/12,-psi1/6,-psi0])
sage: X4 = Z.ambient_space()
sage: Z2 = X4.anticanonical_hypersurface( coefficients=[B/24,1/12,1/3,1/2,B/24,B/12,-psi1/6,-psi0])
sage: Z1 = X4.anticanonical_hypersurface( coefficients=[B/24,1/12,1/3,1/2,B/24,B/12,0,-psi0])
\end{sagecommandline}
The ray matrix of the fan of $X^4$ is now
\begin{equation*}
\sage{X4.fan().rays()}
\end{equation*}
and the fibration $\tld{\beta}\colon X^4 \to B^1$ takes the form
\begin{gather*}
z \mapsto [s : t] = \left[z_0^{12} z_{170} : z_3^{12} z_{168}\right],
\end{gather*}
so working in affine charts we may treat, say, $z_{170}$ as the K3 fibre parameter.

Since the singularities of $Z_2$ are located in charts involving $z_{334}$, it is natural to represent its defining polynomial as
\begin{gather}
h_2 = q_2 z_{334} + r_2,
\label{eq:h2}
\end{gather}
where
\begin{align}
q_2 &=
\frac{1}{12} z_{2}^{12} z_{168}^{11} z_{170}^{11} z_{251}^{8} z_{276}^{4} z_{325}^{6} z_{334}
- \frac{\psi_1}{6} z_{0}^{6} z_{2}^{6} z_{3}^{6} z_{16}^{6} z_{168}^{6} z_{170}^{6} z_{251}^{4} z_{276}^{2} z_{325}^{3}
\nonumber\\
&- \psi_{0} z_{0} z_{1} z_{2} z_{3} z_{4} z_{16} z_{168} z_{170}  z_{251} z_{276} z_{325}
+ \frac{1}{3} z_{1}^{3} z_{251} z_{276}^{2}
+ \frac{1}{2} z_{4}^{2} z_{325},
\label{eq:q2}\\
r_2 &=
\frac{B}{24} z_{16}^{12}
\left(z_3^{12} z_{168} + z_0^{12} z_{170}\right)^{2}.
\label{eq:r2}
\end{align}
In this form it is easy to see that the fibre of $\tld{\beta}\colon Z_2 \to B^1$ (which is generically a K3 surface) over $[s : t] = [-1 : 1]$ splits into two components, corresponding to $q_2 = 0$ and $z_{334} = 0$. The intersection of these components is a curve $C_2$, which is the singular locus of $Z_2$. In the affine chart $(z_1, z_4, z_{170}, z_{334})$ the defining equations of $C_2$ take the form
\begin{align*}
\frac{z_1^3}{3} + \psi_0 z_1 z_4 + \frac{z_4^2}{2} - \frac{\psi_1}{6} = 0,
\qquad
z_{170} = -1,
\qquad
z_{334} = 0,
\end{align*}
which is (generically) a smooth elliptic curve. To see what type of singularities occur in this locus, we translate variables from $(z_1, z_4, z_{170}, z_{334})$ to $(R+u_4, S+u_1, u_2-1, u_3)$, where $(z_1,z_4) = (R, S)$ satisfies the first defining equation of $C_2$ above, to get a hypersurface that is singular at the origin. A further power series substitution brings the singularities into the standard form $u_1^2 + u_2^2 + u_3^2 + (\dots) u_4$, which is a $cA_1$ compound Du Val singularity.  The hypersurfaces $Z_2$ thus have $cA_1$ singularities along $C_2$.

If we now pass to the subfamily $Z_1$, the structure of the singularities remains mostly the same, except that the curve of singularities of $Z_1$, let's call it $C_1$ for this subfamily, develops a singularity of its own and becomes a nodal elliptic curve. In the same affine chart as before the position of the node is $(0,0,-1,0)$. Summarizing, we find:

\begin{proposition} \label{prop:Z1Z2properties} The subfamily of hypersurface models $Z_2$ with $\psi_s = -B$ are generically singular Calabi-Yau threefolds with a smooth elliptic curve $C_2$ of $cA_1$ singularities.

The subfamily of hypersurface models $Z_1$ with $\psi_s = -B$ and $\psi_1 = 0$ are generically singular Calabi-Yau threefolds with a nodal elliptic curve $C_1$ of $cA_1$ singularities.
\end{proposition}

\subsection{Desingularization}
\label{subsection:desingularization}

In this subsection we complete our analysis of the singularities of subfamilies by determining the singularities of the subfamily $Y_1$. We then keep an old promise and use the fibration $\tld{\Phi}\colon X^5_\mpart \to X^4$ to pull back families of hypersurfaces $Z$ to $X^5_\mpart$ then compare them with complete intersections $Y$. 

Since we have changed the resolution used for $X^4$, we need to reconstruct both $X^5_\mpart$ and $\tld{\Phi}$. We do this using $\Sigma^5_\mpart$ as a starting point, constructed in \ref{subsection:ambient space morphism}.
\begin{sagecommandline}
sage: Phi = FanMorphism(m, Sigma5_part, X4.fan(), subdivide=True)
sage: [Delta5_polar.points_pc().index(ray) for ray in Phi.domain_fan().rays()]
[2, 3, 4, 5, 6, 7, 8, 9, 109, 469, 630, 32, 667, 752]
\end{sagecommandline}
We see that $X^5_\mpart$ can still be realized as an open subset of a crepant partial desingularization of $X_{\Delta^5}$. In homogeneous coordinates
\begin{alignat*}{13}
\tld{\Phi}\colon y
&\mapsto
&\big[& z_{0} &:& z_{1} &:& z_{2} &:& z_{3} &:& z_{4} &:& z_{16} &:& z_{168} &:& z_{170} &:& \ \, z_{334} &:& z_{251} &:& z_{276} &:& z_{325}\big]\\
&=
&\big[& y_{4} &:& y_{8} &:& y_{7} &:& y_{5} &:& y_{9} &:& \: y_{6} &:& y_{109} &:& \: y_{32} &:& y_{3}^{2} y_{752} &:& y_{469} &:& y_{630} &:& y_{667}\big].
\end{alignat*}

Considering the subfamily of complete intersections $Y_1 \subset X^5_\mpart$, it is possible to determine that there is a singular point $(0,0,0,0; -1)$ in the chart $(y_2, y_3, y_8, y_9; y_{32})$ (recall that the generating cones of $\Sigma^5_\mpart$ are 4-dimensional, so the corresponding affine charts have a torus factor without a canonical choice of coordinates; in such cases we will use the coordinate corresponding to some suitable ray of the total fan and separate it from the canonically chosen coordinates by ``;''). Note that this point is mapped by $\tld{\Phi}$ to the singular point of the curve of singularities $C_1 \subset Z_1$. It is also a singular point of $X^5_\mpart$ itself, since the cone on rays corresponding to $y_2$, $y_3$, $y_8$, and $y_9$ is not smooth. To fix this we will perform one last subdivision, by inserting the ray corresponding to the midpoint $(-1,-1,1,0,0)$ between ray generators of $y_2$ and $y_3$.
\begin{sagecommandline}
sage: Sigma = Phi.domain_fan().subdivide([(-1,-1,1,0,0)])
sage: all(sigma.is_smooth() for sigma in Sigma if (-1,-1,1,0,0) in sigma)
True
sage: Phi = FanMorphism(m, Sigma, X4.fan())
sage: Phi.is_fibration()
True
sage: points5 = [Delta5_polar.points_pc().index(ray) for ray in Phi.domain_fan().rays()]
sage: points5
[2, 3, 4, 5, 6, 7, 8, 9, 109, 469, 630, 32, 667, 752, 745]
\end{sagecommandline}
The ray matrix of the new fan is
\begin{equation*}
\sage{Phi.domain_fan().rays()}
\end{equation*}
and $\tld{\Phi}$ is still a fibration with coordinate representation
\begin{alignat}{13}
\tld{\Phi}\colon y
&\mapsto
&\big[& z_{0} &:& z_{1} &:& z_{2} &:& z_{3} &:& z_{4} &:& z_{16} &:& z_{168} &:& z_{170} &:& \quad\ z_{334} &:& z_{251} &:& z_{276} &:& z_{325}\big]\nonumber\\
&=
&\big[& y_{4} &:& y_{8} &:& y_{7} &:& y_{5} &:& y_{9} &:& \: y_{6} &:& y_{109} &:& \: y_{32} &:& y_{3}^{2} y_{745} y_{752} &:& y_{469} &:& y_{630} &:& y_{667}\big].
\label{eq:Phi}
\end{alignat}

To construct $Y$ as a subvariety of $X^5_\mpart$ in Sage, we first construct $Y$ in the full space $X^5$ using the same coordinates as for $X^5_\mpart$ (plus $y_0$ and $y_1$), then we use its equations to obtain a subvariety of $X^5_\mpart$. (In the code we refer to it as \sagecommand{Y_part}, but in the text we continue using $Y$ only since mathematically these are the same varieties.)
\begin{sagecommandline}
sage: X5 = CPRFanoToricVariety(np.Delta(), coordinate_points=[0,1]+points5, coordinate_names="y+")
sage: Y = X5.nef_complete_intersection(np, monomial_points= "vertices+origin",coefficients=[[1,1,1],[1,1,1,"xi0","xi1",1]])
sage: X5_part = ToricVariety(Phi.domain_fan(), coordinate_names= "y+", coordinate_indices=points5, base_field=Y.base_ring())
sage: X5_part.inject_coefficients();
sage: X5_part.inject_variables();
sage: S = X5_part.coordinate_ring()
sage: Y_part = X5_part.subscheme([sum(S(g.monomial_coefficient(m)) * S(m.subs(y0=1, y1=1)) for m in g.monomials()) for g in Y.defining_polynomials()])
\end{sagecommandline}
The defining polynomials of $Y$ in the new $X^5_\mpart$ are
\begin{align}
g_0 &= y_5^{12} y_{109} + y_4^{12} y_{32} + y_2 y_3 y_{745},
\label{eq:g0 final}\displaybreak[0]\\
g_1 &=
y_3^2 y_7^{12} y_{109}^{11} y_{469}^8 y_{630}^4 y_{32}^{11} y_{667}^6 y_{752}^2 y_{745}\nonumber\\
&+ \xi_1 y_4^6 y_5^6 y_6^6 y_7^6 y_{109}^6 y_{469}^4 y_{630}^2 y_{32}^6 y_{667}^3 y_{752}
+ \xi_0 y_2^2 y_6^{12} y_{745}\label{eq:g1 final}\\
&+ y_4 y_5 y_6 y_7 y_8 y_9 y_{109} y_{469} y_{630} y_{32} y_{667} y_{752}
+ y_8^3 y_{469} y_{630}^2 y_{752}
+ y_9^2 y_{667} y_{752}.\nonumber
\end{align}

There are \emph{two} singular points on the complete intersections $Y_1$, occurring at the points $\left(\pm 1/\sqrt{\xi_0}, 0, 0, 0; -1\right)$ in the chart $(y_2, y_8, y_9, y_{745}; y_{32})$, which is smooth and does not induce any singularities on its subvarieties. To find their types, consider the defining polynomials of $Y_1$ in this chart:
\begin{align*}
g_0 &=
y_2 y_{745} + y_{32} + 1,\\
g_1 &=
y_{32}^{11} y_{745} + \xi_0 y_2^2 y_{745} + y_{32} y_8 y_9 + y_8^3 + y_9^2,
\end{align*}
then use $g_0$ (which is smooth) to eliminate $y_{32}$ from $g_1$ and shift variables to one of the indicated singular points. The leading terms of the remaining polynomial are quadratic and can be converted into the form $u_1^2 + u_2^2 + u_3^2 + u_4^2$. Both points are therefore conifold singularities on $Y_1$. Thus, in summary, we have:

\begin{proposition} \label{prop:Y1nodal} The subfamily of complete intersection models $Y_1$ with $\xi_1=0$ are generically singular Calabi-Yau threefolds having two isolated nodes.
\end{proposition}

We finally turn our attention to the explicit comparison of $Y$ and the pullback of $Z$ under the fibration $\tld{\Phi}$, combined with the first defining equation of $Y$ (recall that the fibration map $\alpha\colon X^5 \to B^2$ induces the K3 fibration of $Y$ over the curve in $B^2$ corresponding to $g_0 = 0$). Since in \ref{subsection:K3 fibrations} we matched modular parameters of the K3 fibres of $Y$ and $Z_2$, we start with the subfamily $Z_2$.

It is clear that if we use $\tld{\Phi}$~\eqref{eq:Phi} to pullback $h_2$~\eqref{eq:h2} we will not get $g_1$~\eqref{eq:g1 final}, due to a mismatch in the number of monomials and different coefficients, even taking into account the parameter correspondence~\eqref{eq:Y-Z parameter matching}. However, the pullback of $r_2$~\eqref{eq:r2} is
\begin{gather*}
\frac{B}{24} \left(y_5^{12} y_{109} + y_4^{12} y_{32}\right)^2 y_6^{12}
=
\frac{B}{24} \left(y_2 y_3 y_{745}\right)^2 y_6^{12},
\end{gather*}
where equality follows from vanishing of $g_0$~\eqref{eq:g0 final}. Now the number of monomials is no longer an issue, but we can see that $y_3^2 y_{745}$ is a factor of $\tld{\Phi}^*(h_2)$. Further comparison of $\tld{\Phi}^*(h_2)$ and $g_1$ reveals that, subject to $g_0 = 0$,
\begin{gather*}
\Psi^*(h_2) = \frac{y_3^2 y_{745}}{48} \, g_1,
\end{gather*}
where $\Psi$ is $\tld{\Phi}$ precomposed with the following coordinate scaling (which cannot be realized as a toric morphism):
\begin{gather*}
y_6 = \frac{y_6}{\psi_0 \sqrt{12}},
\qquad
y_8 = \frac{y_8}{2},
\qquad
y_9 = - \frac{y_9}{\sqrt{12}},
\qquad
y_{752} = \frac{y_{752}}{2}.
\end{gather*}

Using $\Psi$ to pullback the defining polynomial $h_3$ of a generic member $Z$ of the full 3-parameter family of hypersurfaces, we obtain an isomorphism of $Z$ with the complete intersection defined by $\Psi^*(h_3)$ and $g_0$. Indeed,
\begin{align*}
\Psi^*(h_3)
&=
\Psi^*\left(h_2 - \frac{\psi_s + B}{12} (z_0 z_3 z_{16})^{12} z_{168} z_{170}\right)\\
&=
\frac{y_3^2 y_{745}}{48} \, g_1
- \frac{\psi_s + B}{12^7 \psi_0^{12}} (y_4 y_5 y_6)^{12} y_{32} y_{109},
\end{align*}
and $\Psi^*(h_3) = 0$ implies that $y_3 y_{745} \neq 0$, since otherwise it would be necessary to have $y_4 y_5 y_6 y_{32} y_{109} = 0$, which is impossible: the following computation shows that this set of variables does not appear together in any of the covering charts of $X^5_\mpart$.
\begin{sagecommandline}
sage: charts = [set(S.gen(j) for j in sigma.ambient_ray_indices()) for sigma in X5_part.fan()]
sage: any(chart.intersection([y3, y745]) and chart.intersection([y4, y5, y6, y32, y109]) for chart in charts)
False
\end{sagecommandline}
Thus $g_0 = 0$ can be solved for $y_2$, which is the only fibre variable. This 3-parameter family of complete intersections does not correspond to a nef-partition of $\Delta^5$, the easiest way to see this is to note that it shares one of its equations with $Y$ and, in the case of two-part nef-partitions (and corresponding nef complete intersections), each part completely determines the other.

Coming back to the subfamily $Z_2$, we see that these singular hypersurfaces are pulled back to the union of $Y$ (which is generically smooth) and two toric divisors intersected with $\set{g_0 = 0}$. The composition $\tld{\beta} \circ \Psi$ maps both of these intersections to $[s : t] = [-1 : 1]$. 

Next we look at the preimage under $\Psi$ of the curve of singularities $C_2$ of $Z_2$ inside of $Y$. This preimage satisfies the following conditions:
\begin{gather*}
g_0 = g_1 = \Psi^*(z_{334}) = \Psi^*(r_2) = \Psi^*(q_2) = 0.
\end{gather*}
Conditions $g_0 = \Psi^*(r) = 0$ imply $y_2 y_3 y_{745} = 0$, while $\Psi^*(z_{334}) = 0$ means that $y_3 y_{745} y_{752} = 0$. But $y_2$ and $y_{752}$ cannot vanish simultaneously (we check this below), so we must have $y_3 y_{745} = 0$. We also have $g_1 - 24 y_{752} \Psi^*(q_2) = \xi_0 y_2^2 y_6^{12} y_{745}$ and $y_3$ cannot vanish simultaneously with $y_2$ or $y_6$, so $y_{745} = 0$. Finally, $y_{745}$ and $y_{752}$ cannot vanish simultaneously, so the above conditions are equivalent to 
\begin{gather*}
g_0 = g_1 = y_{745} = 0
\end{gather*}
and we see that the preimage of $C_2$ in $Y$ is a (generically smooth) surface $S_2$. Therefore, we have obtained a geometric transition from $Y$ to $Z$ through the singular subfamily $Z_2$! In fact, $S_2 \simeq C_2 \times \PP^1$ and we have obtained a primitive geometric transition of type III in terms of Definition~\ref{def:primitive transitions}. It just remains to verify the assertions made about the simultaneous vanishing of coordinates:
\begin{sagecommandline}
sage: any(chart.issuperset([y2, y752]) for chart in charts)
False
sage: any(chart.issuperset([y2, y3]) for chart in charts)
False
sage: any(chart.issuperset([y3, y6]) for chart in charts)
False
sage: any(chart.issuperset([y745, y752]) for chart in charts)
False
\end{sagecommandline}

The last thing that we should do in this section is to use the map $\Psi$ to compare the two subfamilies $Z_1$ and $Y_1$. We have already determined that the hypersurfaces $Z_1$ have a curve of singularities $C_1$, which itself has a singular point $(0,0,-1,0)$ in the chart $(z_1, z_4, z_{170}, z_{334})$. Let $S_1 \subset Y_1$ be the preimage of $C_1$ under $\Psi$ (defined by the same equations as $S_2 \subset Y_2$), it is a singular surface. If we consider the preimage of the singular point in $Y_1$, then in addition to the defining equations of $S_1$ we need to impose the condition $y_8 = y_9 = 0$. This leads to
\begin{gather*}
y_8 = y_9 = y_{745} = y_5^{12} y_{109} + y_4^{12} y_{32} = 0,
\end{gather*}
which defines a projective line containing both nodes of $Y_1$. The singular locus of the surface $S_1$ is precisely this line, which is given by $(*,0,0,0;-1)$ in charts $(y_{2 \text{ or } 3}, y_8, y_9, y_{745}; y_{32})$. Thus, we find:

\begin{proposition}\label{prop:geomtrans} The models $Y$ and $Z$ are related by a geometric transition of type III through the singular subfamily $Z_2$. Furthermore, this induces a geometric transition between the subfamilies $Y_1$ and $Z_1$. We have a diagram
\begin{equation}\label{equation:mirror1}\xymatrix{ & Y \ar@{~>}[r] \ar[d] & Y_1 \ar[d] \\
Z \ar@{~>}[r] & Z_2 \ar@{~>}[r] & Z_1}
\end{equation}
where wiggly arrows denote degeneration and straight arrows denote blow-ups.
\end{proposition}

\section{Involutions} \label{section:involutions}

Now we have constructed two families of threefolds, $Y$ and $Z$, along with their singular subfamilies $Y_1$ and $Z_1$, and showed that they are in fact two halves of the same picture: one can move from one to the other by means of a geometric transition. Furthermore, from Proposition \ref{prop:Yproperties} we know that the GKZ series of the singular threefolds with two nodes in the subfamily $Y_1$ agrees with the expected GKZ series for examples of the 14th case. However, we find:

\begin{lemma} \label{lemma:Y1MHS}  The mixed Hodge structure on the cohomology groups $H^i(Y_1)$ of $Y_1$ is pure unless $i = 3$ and the weight filtration $W_{\bullet}$ on $H^3(Y_1)$ has the following description:
\begin{align*} \mathrm{Gr}^{W_{\bullet}}_kH^3(Y_1) &= 0 \quad \mathrm{for}\ k \neq 2,3\\
\dim_{\CC}W_2(H^3(Y_1)) & \in \{1,2\}
\end{align*}
Furthermore, $Y_1$ is factorial if and only if $\dim_{\CC}W_2(H^3(Y_1)) = 2$.
\end{lemma}

\begin{proof} With the exception of the statement about $\dim_{\CC}W_2(H^3(Y_1))$, the statements about the mixed Hodge structure on $H^i(Y_1)$ are an easy consequence of Proposition 3.10 and Corollary 3.13 from \cite{NamikawaSteenbrink:1995}.

To obtain the statement about $\dim_{\CC}W_2(H^3(Y_1))$ we use \cite[Proposition 3.10]{NamikawaSteenbrink:1995}. Recall from Proposition \ref{prop:Y1nodal} that $Y_1$ is singular at precisely two nodes. Let $U$ be a contractible Stein open neighbourhood of either node. Let $\mathrm{Weil}(U)$ (resp. $\mathrm{Cart}(U)$) denote the group of Weil (resp. Cartier) divisors on $U$. Then \cite{Milnor:1968} shows that $\mathrm{Weil}(U)/\mathrm{Cart}(U) \cong \ZZ$. \cite[Proposition 3.10]{NamikawaSteenbrink:1995} thus gives
\[\dim_{\CC}W_2(H^3(Y_1)) = 2 - \sigma(Y_1),\]
where $\sigma(Y_1)$ is the rank of the group of Weil modulo Cartier divisors on $Y_1$. To prove that this is either $1$ or $2$, we simply need to show that $\sigma(Y_1) \neq 2$.

Note that $Y_1$ is smoothable, as it smooths to $Y$. Let $\tilde{Y}_1$ be a small (not necessarily projective) resolution of $Y_1$ with exceptional curves $C_1$ and $C_2$ Then, by \cite[Theorem 2.5]{Namikawa:2002}, there is a non-trivial relation $\alpha_1[C_1] + \alpha_2[C_2] = 0$ in $H_2(\tilde{Y}_1,\CC)$. This implies that $\sigma(Y_1) \neq 2$, as required.

We conclude by showing that $Y_1$ is factorial if and only if $\dim_{\CC}W_2(H^3(Y_1)) = 2$. Note that $\dim_{\CC}W_2(H^3(Y_1)) = 2$ if and only if $\sigma(Y_1) = 0$, which occurs if and only if $Y_1$ is $\QQ$-factorial. But, since $\mathrm{Weil}(U)/\mathrm{Cart}(U)$ is torsion-free in a neighbourhood of each node, $Y_1$ is $\QQ$-factorial if and only if it is factorial.\end{proof}

In order to show which of the two cases from Lemma \ref{lemma:Y1MHS} hold, we will prove that $Y_1$ is factorial. In order to do this, we note that it is easier to show factoriality for the single node family of threefolds obtained as the quotient of $Y_1$ by a certain involution. We will then use factoriality of this family to deduce factoriality of $Y_1$. To this end, we begin this section by studying several involutions on the complete intersections $Y$. 

One of these involutions, $i_1$, is easily visible from the vertices of $\left(\Delta^5\right)^\circ$ or rays of $\Sigma^5$: it is the exchange of the first two affine coordinates. In homogeneous coordinates it is realized by
\begin{gather*}
y_4 \leftrightarrow y_5,
\qquad
y_{32} \leftrightarrow y_{109},
\end{gather*}
with the other $y_i$ unchanged, and in terms of the K3 fibration of $Y$ over $\PP^1$ induced by $\tld{\alpha}$, as discussed in \ref{subsection:K3 fibrations}, it corresponds to $[u : v] \leftrightarrow [v : u]$.

Another involution, $i_2$, is  given in homogeneous coordinates \emph{on the unresolved toric variety $X_{\Delta^5}$} by
\begin{gather*}
y_2 \leftrightarrow y_3,
\qquad
y_6 \leftrightarrow y_7,
\end{gather*}
with the other $y_i$ unchanged, followed by the rescaling map
\begin{gather*}
y_2 \mapsto \frac{y_2}{\sqrt{\xi_0}},
\qquad
y_3 \mapsto y_3 \sqrt{\xi_0}.
\end{gather*}
The resolution that we are currently using for $X^5$ is not suitable for a clean realization of the first part of $i_2$ as a toric morphism, however for a suitable refinement of the fan the underlying lattice map is given by the matrix
\begin{sagesilent}
m_1 = Delta5_polar.vertices_pc()(2,3,4,5,6,7,8,9).matrix()
m_2 = Delta5_polar.vertices_pc()(3,2,4,5,7,6,8,9).matrix()
mN = m_2.solve_right(m_1)
\end{sagesilent}
\begin{equation} \label{eq:mN}
m_N = \sage{mN}.
\end{equation}
For ``unsuitable'' refinements we can view $i_2$ as a birational map, which is certainly defined on the torus and may extend to some of the lower dimensional orbits as well. This involution corresponds to the same involution on the $\PP^1$ base of one of the elliptic fibrations on the $M$-polarized K3 fibres as $i_1$ did on the $\PP^1$ base of the K3 fibration of $Y$.

Yet another involution, $i_3$, is given in homogeneous coordinates by
\begin{gather*}
y_9
\mapsto
- y_9
- y_4 y_5 y_6 y_7 y_8 y_{109} y_{469} y_{630} y_{32},
\end{gather*}
with the other $y_i$ unchanged. It comes from either of the two elliptic fibrations on the $M$-polarized K3 surfaces, as the $y \mapsto -y$ involution in the coordinates of the Weierstrass normal form.  This is  the same involution as the one described in \cite[Section~4.2.1]{Billo:1998}.\footnote{There it is realized as $x_5 \mapsto - x_5$ in the ``alternative gauge''. One can convert it to the gauge used in \ref{subsection:AH model} via formulas (4.13) and (4.15) in~\cite{Billo:1998} (note that those formulas contain a typo in the expression for $\tld{\psi}_4$: the term $- 4 b_4 {\lambda'_4}^2 \lambda_4 \psi_0$ should not contain $\psi_0$, while the term $4 \lambda'_5 \lambda'_4 \lambda_{03}$ should).}

It follows from the representations of these three involutions in homogeneous coordinates that they commute with each other, thus it is possible to take their compositions to obtain new involutions, forming a group isomorphic to $\ZZ_2^3$.

\subsection{Action on the Singular Locus}

Of particular interest is the action of the involutions on the singular locus of $Y_1$ and its related subvarieties. To analyse this action we can work in the chart $(y_2, y_8, y_9, y_{745}; y_{32})$, used before. 

It is very easy to restrict the last involution, $i_3$, to the chosen chart; we obtain $y_9 \mapsto - y_9 - y_8 y_{32}$. On the exceptional K3 fibre $y_{32} = - 1$ containing the nodes, this becomes $y_9 \mapsto y_8 - y_9$. The line passing through the nodes with $y_8 = y_9 = 0$ is a part of the fixed point locus. 

In order to get the chart representation of $i_1$ we need to know the relationship between $y_{109}$ and the variables of the chosen chart.
\begin{sagecommandline}
sage: Delta5_polar.points_pc()(2,8,9,745,32,109).matrix().kernel()
Free module of degree 6 and rank 1 over Integer Ring
Echelon basis matrix:
[ 11   4   6 -10   1   1]
\end{sagecommandline}
This shows that there is a scaling action by $\CC^*$ on the homogeneous coordinates $y$ given by
\begin{gather*}
\left[ \lambda^{11} y_2 : \lambda^4 y_8 : \lambda^6 y_9 : \lambda^{-10} y_{745} : \lambda y_{32} : \lambda y_{109} \right].
\end{gather*}
Using this to eliminate $y_{109}$, we get 
\begin{gather*}
\left[ y_{109}^{-11} y_2 : y_{109}^{-4} y_8 : y_{109}^{-6} y_9 : y_{109}^{10} y_{745} : y_{109}^{-1} y_{32} : 1 \right],
\end{gather*}
so the chart representation of $i_1$ is
\begin{gather*}
(y_2, y_8, y_9, y_{745}; y_{32})
\mapsto
\left( y_{32}^{-11} y_2, y_{32}^{-4} y_8, y_{32}^{-6} y_9, y_{32}^{10} y_{745}; y_{32}^{-1} \right),
\end{gather*}
which reduces to
\begin{gather*}
(y_2, y_8, y_9, y_{745}; -1)
\mapsto
\left( - y_2, y_8, y_9, y_{745}; -1 \right),
\text{ i.e. }
y_2 \mapsto - y_2,
\end{gather*}
on the exceptional K3 fibre containing the singular locus. Note that this action exchanges the nodes of $Y_1$.

We now perform a similar computation for $i_2$ (the first step is to determine the indices of points of $\left(\Delta^5\right)^\circ$ which are images under the action by the matrix~\eqref{eq:mN} of the points corresponding to our current homogeneous coordinates $y$):
\begin{sagecommandline}
sage: points5
sage: [Delta5_polar.points_pc().index(Delta5_polar.point(p)*mN) for p in points5]
sage: Delta5_polar.points_pc()(2,8,9,745,32,3).matrix().kernel()
Free module of degree 6 and rank 1 over Integer Ring
Echelon basis matrix:
[ 1  0  0 -2  0  1]
sage: Delta5_polar.points_pc()(2,8,9,745,32,21).matrix().kernel()
Free module of degree 6 and rank 1 over Integer Ring
Echelon basis matrix:
[ 11  0  0 -11  1  -1]
\end{sagecommandline}
leading to the chart representation
\begin{gather*}
i_2\colon (y_2, y_8, y_9, y_{745}; y_{32})
\mapsto
\left( \xi_0^{-1} y_2^{-1} y_{32}^{11}, y_8, y_9, \xi_0 y_2^2 y_{32}^{-11} y_{745}; y_{32} \right).
\end{gather*}
On the exceptional K3 fibre this reduces to
\begin{gather*}
i_2\colon (y_2, y_8, y_9, y_{745}; -1)
\mapsto
\left( - \xi_0^{-1} y_2^{-1}, y_8, y_9, - \xi_0 y_2^2 y_{745}; -1 \right)
\end{gather*}
which also exchanges the nodes of $Y_1$.

\subsection{Fixed Point Loci} \label{subsection:fixed loci}

We will determine the fixed point loci of the involutions in the same affine chart as before. The defining polynomials of $Y$ in this chart are
\begin{align*}
g_0 &=
y_2 y_{745} + y_{32} + 1,\\
g_1 &=
y_{32}^{11} y_{745} + \xi_1 y_{32}^6 + \xi_0 y_2^2 y_{745} + y_{32} y_8 y_9 + y_8^3 + y_9^2.
\end{align*}

For $i_1$ it is easy to see that the fixed point locus has two components, one is determined by $y_{32} = 1$ and the other by $y_{32} = -1, y_2 = 0$; both correspond to surfaces on $Y$.

For $i_2$ the fixed point locus is determined by
\begin{gather*}
y_2 = \xi_0^{-1} y_2^{-1} y_{32}^{11}
\qquad\text{and}\qquad
y_{745} = \xi_0 y_2^2 y_{32}^{-11} y_{745},
\end{gather*}
which is equivalent to the single condition
\begin{gather*}
\xi_0 y_2^2 = y_{32}^{11},
\end{gather*}
which cuts out a surface in $Y$.

Finally, for $i_3$ the fixed point locus is defined by $2 y_9 + y_8 y_{32} = 0$, which again determines a surface in $Y$.

Each of the $i_j$ changes the sign of the holomorphic $(3,0)$-form on $Y$ (this can be seen from \cite[Section 6]{Billo:1998}), so an interesting subgroup of the involutions is formed by their pairwise compositions, which preserve the sign of the holomorphic $(3,0)$-form.

The composition $i_2 \circ i_3$ gives, in fact, a Nikulin involution on the smooth K3 fibres of $Y$ with exactly $8$ fixed points in each of them. More explicitly,
\begin{gather*}
i_2 \circ i_3\colon  (y_2, y_8, y_9, y_{745}; y_{32})
\mapsto
\left( \xi_0^{-1} y_2^{-1} y_{32}^{11}, y_8, - y_9 - y_8 y_{32}, \xi_0 y_2^2 y_{32}^{-11} y_{745}; y_{32} \right),
\end{gather*}
so the fixed point locus is determined by intersection of the fixed point loci of $i_2$ and $i_3$ considered separately:
\begin{gather*}
\xi_0 y_2^2 = y_{32}^{11}
\qquad\text{and}\qquad
2 y_9 + y_8 y_{32} = 0.
\end{gather*}
Specifying the K3 fibration base coordinate $y_{32}$ generically gives a two-fold ambiguity in $y_2$, which allows us to uniquely determine $y_{745}$ from $g_0$ and get a cubic equation for $y_8$ from $g_1$. This gives us $6$ fixed points: two more should be located outside this chart. It does not seem possible, however, to locate them in the other charts of our resolution of $X_{\Delta^5}$, due to the fact that $i_2$ is only a birational map, not a morphism.

Let us work out the 6 visible points in more detail treating $y_{32} \in \CC^*$ as a parameter; this will be useful for the analysis we will perform later in \ref{subsection:construction}. In this case
\begin{align*}
y_2 &= \pm \frac{y_{32}^6}{\sqrt{\xi_0 y_{32}}},\\
y_{745} &= - \frac{y_{32} + 1}{y_2}, \\
0 &= y_8^3 - \frac{y_{32}^2}{4} y_8^2 + 2 y_{32}^6 \left(\frac{\xi_1}{2}  \mp \sqrt{\xi_0 y_{32}} \mp \frac{\xi_0}{\sqrt{\xi_0 y_{32}}}\right), \\
y_9 &= - \frac{y_8 y_{32}}{2},
\end{align*}
with a consistent choice of signs. 

Next we consider the composition
\begin{gather*}
i_1 \circ i_2\colon  (y_2, y_8, y_9, y_{745}; y_{32})
\mapsto
\left( \xi_0^{-1} y_2^{-1}, y_{32}^{-4} y_8, y_{32}^{-6} y_9, \xi_0 y_2^2 y_{32}^{-1} y_{745}; y_{32}^{-1} \right).
\end{gather*}
Its fixed point locus is concentrated in fibres corresponding to $y_{32} = \pm 1$, here the involution takes the form
\begin{gather*}
(y_2, y_8, y_9, y_{745}; \pm 1)
\mapsto
\left( \xi_0^{-1} y_2^{-1}, y_8, y_9, \pm \xi_0 y_2^2 y_{745}; \pm 1 \right),
\end{gather*}
where all alternating signs must be chosen to be the same. Thus the fixed points satisfy the conditions
\begin{gather*}
\xi_0 y_2^2 = 1
\qquad\text{and}\qquad
\xi_0 y_2^2 y_{745} = \pm y_{745}.
\end{gather*}
In the case $y_{32} = - 1$ the second condition yields $y_{745} = 0$, which also follows from $g_0 = 0$, so on $Y$ the fixed point locus is determined by $y_{32} = \pm 1$ and $\xi_0 y_2^2 = 1$. This corresponds to four elliptic curves, two of which become nodal with nodes coinciding with the nodes of $Y_1$ in the case $\xi_1 = 0$.

The last pairwise composition is
\begin{gather*}
i_1 \circ i_3\colon  (y_2, y_8, y_9, y_{745}; y_{32})
\mapsto
\left( y_{32}^{-11} y_2, y_{32}^{-4} y_8, - y_{32}^{-6} y_9 - y_{32}^{-5} y_8, y_{32}^{10} y_{745}; y_{32}^{-1} \right).
\end{gather*}
Its fixed point locus is also located in K3 fibres with $y_{32} = \pm 1$. For $y_{32} = 1$ the only other condition is $2 y_9 + y_8 = 0$, while for $y_{32} = - 1$ we need both $2 y_9 - y_8 = 0$ and $y_2 = 0$. Taking into account the defining equations of $Y$, we see that we for $y_{32} = 1$ we get an elliptic curve, while for $y_{32} = -1$ the component of the fixed point locus is a 3-section of the elliptic fibration on the surface defined by $y_2 = 0$.

We also observe that the intersection of the fixed point loci of these compositions on $Y$ is given by
\begin{gather*}
y_{32} = 1,
\qquad
\xi_0 y_2^2 = 1,
\qquad
2 y_9 + y_8 = 0,
\end{gather*}
which coincides with the fixed point locus of the Nikulin involution in the K3 fibre corresponding to $y_{32} = 1$.

\subsection{Quotient by the Nikulin involution}

By the results of \ref{subsection:fixed loci}, a complete intersection $Y$ admits a fibrewise Nikulin involution $i_2 \circ i_3$. If we quotient $Y_1$ by this involution and perform a crepant partial desingularization of the resulting quotient singularities, we obtain a new Calabi-Yau threefold $W$. By \cite[Theorem 3.13]{ClingherDoran:2007}, the K3 fibration on $Y_1$ induces a fibration on $W$ whose general fibres are Kummer surfaces and, furthermore, if $\{E_1,E_2\}$ are the pair of elliptic curves canonically associated to a fibre of $Y_1$ by Theorem \ref{theorem:M-polarized K3}, the corresponding fibre of $W$ is $\mathrm{Kum}(E_1 \times E_2)$ (i.e. the Kummer surface associated to $E_1 \times E_2$). This threefold $W$ will be constructed in a more direct way in \ref{section:forward construction}.

\begin{proposition} \label{prop:WMHS} $W$ is a Calabi-Yau threefold that moves in a $1$-parameter family, the general member of which is factorial and smooth away from a single node. The mixed Hodge structure on the cohomology groups $H^i(W)$ of $W$ is pure unless $i = 3$ and the weight filtration $W_{\bullet}$ on $H^3(W)$ has the following description:
\begin{align*} \mathrm{Gr}^{W_{\bullet}}_kH^3(W) &= 0 \quad \mathrm{for}\ k \neq 2,3\\
\dim_{\CC}W_2(H^3(W)) & = 1
\end{align*}
\end{proposition}
\begin{proof} Note that quotient of $Y_1$ by the fibrewise Nikulin involution introduces  a curve of $A_1$ singularities, which are resolved by a single crepant blow-up. The two nodes of $Y_1$ are identified by the quotient to give a single node on $W$.

With this in place, the remainder of the proof proceeds in the same way as the proof of Lemma \ref{lemma:Y1MHS}. Note here that the involution $i_2 \circ i_3$ is defined on the family of smooth threefolds $Y$, not just on $Y_1$, and the resolved quotient of $Y$ by this involution gives a smoothing of $W$.
\end{proof}

It is immediate from factoriality that $W$ does not admit a small projective resolution. In fact, more is true:

\begin{corollary} \label{corollary:Wsympres} $W$ does not admit a small symplectic resolution. In particular, the singularity of $W$ cannot be resolved to give a smooth Calabi-Yau threefold.
\end{corollary}
\begin{proof} Since $W$ does not admit a small projective resolution, \cite[Corollary 8]{CortiSmith:2005} shows that the Lagrangian $3$-sphere obtained by smoothing the single node cannot be nullhomologous. But, by \cite[Theorem 2.9]{STY:2002}, this implies that $W$ does not admit a small symplectic resolution.

The statement about non-existence of Calabi-Yau resolutions is immediate from the fact that $W$ has terminal singularities and no small resolution. Thus any resolution of $W$ will contribute positively to the canonical divisor, violating the Calabi-Yau condition.
\end{proof}

From Proposition \ref{prop:WMHS} we can deduce factoriality of $Y_1$. This will in turn allow us to show that $Y_1$ provides the geometric realisation of the 14th case variation of Hodge structure that we have been searching for.

\begin{theorem} \label{thm:14thcase} The singular Calabi-Yau threefold $Y_1$ is factorial, so the weight filtration $W_{\bullet}$ on $H^3(Y_1)$ has $\dim_{\CC}W_2(H^3(Y_1)) = 2$. Furthermore, as $Y_1$ varies in its one-parameter subfamily, the graded piece $\mathrm{Gr}^{W_{\bullet}}_3H^3(Y_1)$ admits a pure variation of Hodge structure of weight $3$ and type $(1,1,1,1)$, which realizes the 14th case variation of Hodge structure.
\end{theorem}

\begin{proof} Clingher and Doran \cite{ClingherDoran:2007}\cite{ClingherDoran:2011} have shown that the fibrewise Nikulin involution $i_2 \circ i_3$ is a \emph{geometric $2$-isogeny} on the fibres of $Y_1$, i.e. there exists a second fibrewise Nikulin involution on $W$ such that the resolved quotient of $W$ by this second involution is isomorphic to $Y_1$. As $W$ is factorial and $\QQ$-factoriality is preserved under finite surjective morphisms \cite[Lemma 5.16]{KollarMori:1998} and resolutions, we thus have that $Y_1$ must be $\QQ$-factorial. Then, from the proof of Lemma  \ref{lemma:Y1MHS}, we see that $Y_1$ is $\QQ$-factorial if and only if it is factorial.

Given this, the fact that $\dim_{\CC}W_2(H^3(Y_1)) = 2$ follows from Lemma \ref{lemma:Y1MHS}. The statement about the Hodge decomposition of $\mathrm{Gr}^{W_{\bullet}}_3H^3(Y_1)$ then follows from calculation of $b_3(Y_1) = 6$, using \cite[Theorem 3.2] {NamikawaSteenbrink:1995}.

Finally, recall from Proposition \ref{prop:Yproperties} that the GKZ series of $Y_1$ matches that predicted for the 14th case VHS. Furthermore, it follows from the classification of Doran and Morgan \cite{DoranMorgan:2006} that the 14th case VHS is the only one with $h^{2,1} = 1$ having this GKZ series. The variation of Hodge structure on $\mathrm{Gr}^{W_{\bullet}}_3H^3(Y_1)$ must therefore be the 14th case VHS.
\end{proof}

Again, factoriality implies that $Y$ does not admit a small projective resolution. Indeed, as for $W$, a stronger statement can be made:

\begin{corollary} \label{corollary:Y1sympres} $Y_1$ does not admit a small symplectic resolution. In particular, the singularities of $Y_1$ cannot be resolved to give a smooth Calabi-Yau threefold.
\end{corollary}
\begin{proof} Given that $b_3(Y)=8$ and $b_3(Y_1) = 6$, \cite[Theorem 3.2]{Rossi:2006} shows that the two vanishing cycles in $Y$ are homologically independent. But, by \cite[Theorem 2.9]{STY:2002}, this implies that $Y_1$ does not admit a small symplectic resolution. The statement about non-existence of Calabi-Yau resolutions is proved in exactly the same way as the corresponding statement in Corollary \ref{corollary:Wsympres}.
\end{proof}
\begin{remark} We note that the statement of \cite[Theorem 3.2]{Rossi:2006} technically requires that $Y_1$ have a small projective resolution. However, by examining the proof it is easy to see that the assertion about the number of homologically independent vanishing cycles is independent of this assumption.
\end{remark}

\section{The Forward Construction} \label{section:forward construction}
\subsection{Undoing the Kummer construction.} \label{subsection:construction}

The one-parameter family of Calabi-Yau threefolds $W$ is very closely related to $Y_1$; for instance, the quotient-resolution procedure does not affect the GKZ series, so $W$ provides a second example of a one-parameter family with GKZ series equal to that predicted for the 14th case. It could thus prove fruitful to study the geometry of $W$ explicitly. In order to do this, we apply the methods of \cite{DHNT:2013} to produce an explicit model for $W$ and, in doing so, illustrate the usage of these methods in a concrete example.

Recall that $W$ admits a fibration whose general fibres are Kummer surfaces and, furthermore, if $\{E_1,E_2\}$ are the pair of elliptic curves canonically associated to a fibre of $Y_1$ by Theorem \ref{theorem:M-polarized K3}, the corresponding fibre of $W$ is $\mathrm{Kum}(E_1 \times E_2)$ (i.e. the Kummer surface associated to $E_1 \times E_2$). Our aim is to find a method by which $W$ can be constructed directly by performing a fibrewise Kummer construction on a threefold fibred by products of elliptic curves.

In order to do this, we begin by using the methods of \cite[Section 4]{DHNT:2013} to undo the Kummer construction on $W$, i.e. find a family of Abelian surfaces $\cA$, such that the application of the Kummer construction fibrewise to $\cA$ yields a birational model of $W$. Unfortunately, as we shall see, this is not possible directly: instead, we will need to proceed to a finite cover of $W$ first.

We begin with some setup. Let $U \subset \PP^1$ be the open set over which the fibres of the K3 fibration $Y_1 \to \PP^1$ induced by $\tilde{\alpha}$ are nonsingular and let $p \in U$ be a general point. Applying the method of \cite[Section 4.3]{Billo:1998} to the affine description of $Y_1$ given in \ref{subsection:fixed loci} (and recalling that $Y_1$ is obtained by setting $\xi_1 = 0$ in the equations for $Y$), we see that the singular fibres of $Y_1$ occur over $y_{32} \in \{0,-1,\alpha,\beta,\infty\}$ where, as before, $y_{32}$ is an affine coordinate on the base $\PP^1$ of the K3 fibration, and $\alpha$ and $\beta$ are given by
\[\alpha,\beta = \frac{2 - 12^6 \xi_0 \pm 2 \sqrt{1 - 12^6 \xi_0}}{12^6 \xi_0}.\]
For now we will assume $\xi_0 \notin \{0, \frac{1}{12^6}, \infty\}$, as degenerate behaviour occurs at these points, although the case $\xi_0 = \frac{1}{12^6}$ will be discussed later in \ref{subsection:degenerate}. $U$ is thus a copy of $\PP^1$ with these five points removed.

Let $Y_{1,U} \to U$ denote the restriction of $Y_1 \to \PP^1$ to $U$. Then we have:

\begin{lemma} $Y_{1,U}$ is an $M$-polarized family of K3 surfaces over $U$, in the sense of \cite[Definition 2.1]{DHNT:2013}.
\end{lemma}
\begin{proof} From the discussion in \ref{subsection:K3 fibrations}, we already know that the fibres of $Y_{1,U}$ are $M$-polarized. In order to complete the proof of the lemma, we just need to show that the polarization is invariant under monodromy around the punctures in $U$. To do this, we first note that the lattice $M$ is generated by the classes of the divisors given in \ref{figure:G}. As these divisors are toric, they are invariant under monodromy in $U$. Thus the lattice $M$ is also.
\end{proof}

Proceeding now to $W$, we see that the restricted family of K3 surfaces $W_U \to U$ has fibres isomorphic to $\mathrm{Kum}(E_1 \times E_2)$. Let $\mathrm{NS}(W_p)$ denote the N\'{e}ron-Severi lattice of the fibre of $W_U \to U$ over $p$. Then \cite[Theorem 3.3]{DHNT:2013} shows that $W_U$ is an $(\mathrm{NS}(W_p),G)$-polarized family of K3 surfaces, in the sense of \cite[Definition 2.4]{DHNT:2013}, where $G$ is a finite group controlling the action of monodromy in $U$ on $\mathrm{NS}(W_p)$.

By \cite[Proposition 4.1]{DHNT:2013}, we may undo the Kummer construction on $W$ if $G$ is trivial. Unfortunately, this will prove not to be the case, so we will need to proceed to a finite cover of $W$ to kill the action of $G$. In order to do this, we use the methods of \cite[Section 4.3]{DHNT:2013} to find a finite cover $V \to U$ so that the action of monodromy in $V$ on $\mathrm{NS}(W_p)$ is trivial.

We begin by noting that, by the discussion in \cite[Section 4.3]{DHNT:2013}, the action of monodromy on $\mathrm{NS}(W_p)$ may be computed from its action on the eight exceptional curves $\{F_1,\ldots,F_8\}$ in $W_p$ arising from the blow-up of the eight fixed points of the Nikulin involution.

The arrangement of these eight curves may be described in terms of a certain elliptic fibration on $W_p$. By \cite[Proposition 3.10]{ClingherDoran:2007}, the fibrewise Nikulin involution $i_2\circ i_3$ on $Y_1$ acts on the alternate fibration as translation by a section, so the alternate fibration induces an elliptic fibration $\Psi$ on $W_p$. As the alternate fibration is invariant under monodromy in $Y_{1,U}$, the induced fibration $\Psi$ is invariant under mondromy in $W_U$. Furthermore, the fact that $\sigma = 1$ for $Y_1$ implies that \cite[Assumption 4.6]{DHNT:2013} holds, so $\Psi$ has four sections and seven singular fibres: one of Kodaira type $I_6^*$ and six of Kodaira type $I_2$.

Given this description, we can locate the exceptional curves $F_i$ in relation to the fibration $\Psi$. The $I_6^*$ fibre contains two such curves, which we label $F_1$ and $F_2$, and each $I_2$ fibre contains precisely one such curve, giving $\{F_3,\ldots,F_8\}$. Furthermore, of the four sections
\begin{itemize}
\item one intersects none of the $F_i$,
\item one intersects all of $\{F_3,\ldots,F_8\}$,
\item one intersects $F_1$ and precisely three of $\{F_3,\ldots,F_8\}$ (say $F_3$, $F_4$, $F_5$), and
\item one intersects $F_2$ and the other three $\{F_3,\ldots,F_8\}$ (say $F_6$, $F_7$, $F_8$).
\end{itemize}

By \cite[Proposition 4.8]{DHNT:2013}, the action of monodromy around a loop in $U$ must either
\begin{enumerate}
\item fix both $F_1$ and $F_2$, in which case the sets $\{F_3,F_4,F_5\}$ and $\{F_6,F_7,F_8\}$ are preserved, or
\item swap $F_1$ and $F_2$, in which case the sets $\{F_3,F_4,F_5\}$ and $\{F_6,F_7,F_8\}$ are interchanged.
\end{enumerate}

In order to see which of these cases occurs, we will compute the action of monodromy on the two sets $\{F_3,F_4,F_5\}$ and $\{F_6,F_7,F_8\}$. For clarity of notation, we will use coordinates $[u:v]$ on the $\PP^1$ base of the fibrations $W \to \PP^1$ and $Y_1 \to \PP^1$, where the affine coordinate $y_{32} = \frac{u}{v}$; note that this agrees with the notation of equations \eqref{eq:pi_Y} and \eqref{eq:sigma_Y}. We have:

\begin{proposition} \label{proposition:Wmonodromy} The action of monodromy in $U$ on $\mathrm{NS}(W_p)$ around:
\begin{itemize}
\item $[u:v] \in \{[-1:1],[\alpha:1],[\beta:1]\}$ fixes both $F_1$ and $F_2$.
\item $[u:v] \in \{[0:1],[1:0]\}$ interchanges $F_1$ and $F_2$.
\end{itemize}
\end{proposition}
\begin{proof} As noted above, it is enough to show that the two sets $\{F_3,F_4,F_5\}$ and $\{F_6,F_7,F_8\}$ are preserved by monodromy around $\{[-1:1],[\alpha:1],[\beta:1]\}$ and interchanged by monodromy around $\{[0:1],[1:0]\}$. Furthermore, recall that the $F_i$ arise as the blow up of the eight fixed points of the Nikulin involution, so in order to understand their behaviour under monodromy on $W_U$ it is sufficient to study the action of monodromy on $Y_{1,U}$ on these eight fixed points.

In \ref{subsection:fixed loci} an affine description for $Y_1$ was computed and the locations of the six fixed points corresponding to the curves $F_3,\ldots,F_8$ within it were calculated explicitly. In particular, we found that the $y_2$-coordinates of these fixed points satisfy
\[y_2 = \pm \frac{y_{32}^6}{\sqrt{\xi_0y_{32}}},\]
with each choice of sign corresponding to three fixed points. These two sets of three fixed points are precisely those corresponding to the two sets of curves $\{F_3,F_4,F_5\}$ and $\{F_6,F_7,F_8\}$. Therefore, these two sets are interchanged by monodromy if and only if monodromy switches the sign of $y_2$ in the expression above. But this happens precisely for monodromy around $y_{32} \in \{0,\infty\}$, i.e. around $[u:v] \in \mbox{\{[0:1],[1:0]\}}$.
\end{proof}

We therefore see that $\mathrm{NS}(W_p)$ is \emph{not} fixed under monodromy in $U$. This presents an obstruction to undoing the Kummer construction on $W$. To resolve this, we pass to a cover of $W$.

Define $W'$ to be the threefold fibred over $\PP^1$ (with base coordinate $[s:t]$) obtained as the pull-back of $W \to \PP^1$ by the map $f \colon  \PP^1 \to \PP^1$ given by $[s:t] \mapsto [s^2:t^2] = [u:v]$. The fibration $W' \to \PP^1$ has eight singular fibres, over 
\[[s:t] \in \{[0:1],[i:1],[-i:1],[\alpha':1],[-\alpha':1],[\beta':1],[-\beta':1],[1:0]\},\]
where $(\alpha')^2 = \alpha$ and $(\beta')^2 = \beta$. Let
\[U' = \PP^1 \setminus  \{[0:1],[i:1],[-i:1],[\alpha':1],[-\alpha':1],[\beta':1],[-\beta':1],[1:0]\};\]
note that $U' = f^{-1}(U)$.

Let $W'_p$ denote a fibre of $W'$ over a point $p \in U'$. Then $W'_p$ is smooth and isomorphic to $W_{f(p)}$, so we may identify the eight curves $\{F_1,\ldots,F_8\}$ in $W'_p$ and, as before, the action of monodromy in $U'$ on $\mathrm{NS}(W'_p)$ may be computed from its action these eight curves. Furthermore, by construction the action of monodromy in $U'$ on $\mathrm{NS}(W'_P)$ fixes $F_1$ and $F_2$ and so, by \cite[Proposition 4.8]{DHNT:2013}, the sets $\{F_3,F_4,F_5\}$ and $\{F_6,F_7,F_8\}$ are preserved. 

It remains to compute the action of monodromy on these two sets. To do this, we use the method from the proof of Proposition \ref{proposition:Wmonodromy} and track the locations of the fixed points of the Nikulin involution. Taking the double cover $W' \to W$ corresponds to making a change of coordinates $y_{32} = \frac{u}{v} = \frac{s^2}{t^2} =: x^2$ in the equations for these fixed points in the affine description of \ref{subsection:fixed loci}, giving
\begin{align*}
y_2 &= \pm \frac{x^{11}}{\sqrt{\xi_0}}, \\
y_{745} &= \mp \sqrt{\xi_0} \frac{x^2 +1}{x^{11}},\\
0 &= y_8^3 - \frac{1}{4}x^4y_8^2 \mp 2\sqrt{\xi_0} x^{11} \mp 2\sqrt{\xi_0} x^{13},\\
y_9 &= - \frac{y_8x^2}{2},
\end{align*}
where we make a consistent choice of signs throughout.

For a given choice of signs, there are three solutions to these equations. These three solutions correspond to the three divisors $\{F_3,F_4,F_5\}$ for one choice of sign and the divisors $\{F_6,F_7,F_8\}$ for the other. Once a choice of signs has been made, $y_2$ and $y_{745}$ are uniquely determined by the first two equations, the third equation gives three possibilities for $y_8$, and the fourth equation uniquely determines $y_9$. Thus we see that the action of monodromy on these equations is completely determined by its action on the third equation. We have therefore shown:

\begin{lemma} The action of monodromy in $U'$ on the set of divisors $\{F_3,F_4,F_5\}$ is identical to its action on the roots of the cubic equation
\[0 = t^{13}y_8^3 - \frac{1}{4}s^4t^9y_8^2  - 2\sqrt{\xi_0} s^{11}t^2 - 2\sqrt{\xi_0} s^{13}\]
and its action on the set of divisors $\{F_6,F_7,F_8\}$ is identical to its action on the roots of the cubic equation
\[0 = t^{13}y_8^3 - \frac{1}{4}s^4t^9y_8^2 + 2\sqrt{\xi_0} s^{11}t^2 + 2\sqrt{\xi_0} s^{13}.\]
\end{lemma}

This monodromy can be calculated using the \texttt{monodromy} command in \emph{Maple}'s \texttt{algcurves} package (this calculation was performed with the generic value $\xi_0 = 1$). Using the base point $[s:t] = [-1:10]$, \ref{table:monodromy} shows the action of monodromy around anticlockwise loops about each of the eight singular fibres in $W'$ on each of the sets $\{F_3,F_4,F_5\}$ and $\{F_6,F_7,F_8\}$, expressed as permutations in $S_3$ (in cycle notation). It is easy to see that the permutations listed in \ref{table:monodromy} generate the full group $S_3 \times S_3$, so monodromy acts as $S_3 \times S_3$ on the sets $\{F_3,F_4,F_5\}$ and $\{F_6,F_7,F_8\}$ in $W'_{[-1:10]}$. 

\begin{table}
\begin{tabular}{|c|c|c|}
\hline 
Monodromy around $[s:t]$ & Action on $\{F_3,F_4,F_5\}$ & Action on  $\{F_6,F_7,F_8\}$ \\
\hline
$[0:1]$ & $(123)$ & $(456)$ \\
$[i:1]$ & $(12)$ & $(45)$ \\
$[-i:1]$ & $(13)$ & $(45)$ \\
 $[\alpha':1]$ & $(23)$ & id \\
$[-\alpha':1]$ & id & $(46)$ \\
$[\beta':1]$ & $(23)$ & id \\
 $[-\beta':1]$ & id & $(56)$ \\
$[1:0]$ & $(123)$ & $(456)$\\
\hline
\end{tabular}
\caption{Action of monodromy on $\{F_3,F_4,F_5\}$ and $\{F_6,F_7,F_8\}$.}
\label{table:monodromy}
\end{table}

Therefore, in order to undo the Kummer construction, we need to pull-back $W'$ to a further $36$-fold cover $g \colon C \to \PP^1$. This cover is constructed as follows: the $36$ preimages of the point $[-1:10] \in \PP^1$ are labelled by elements of $S_3 \times S_3$ and monodromy around each of the points $[s:t]$ listed in \ref{table:monodromy} acts on these labels as composition with the corresponding permutation. This action extends to an action of $S_3 \times S_3$ on the whole of $C$. The map $g$ has ramification index $2$ at all points over $[s:t] \in \{[i:1],[-i:1],[\alpha':1],[-\alpha':1],[\beta:1],[-\beta':1]\}$ and ramification index $3$ at all points over $[s:t] \in \{[1:0],[0:1]\}$.

Let $W'' \to C$ denote the pull-back of $W' \to \PP^1$ under this $36$-fold cover, let $V = g^{-1}(U')$ and let $W''_V$ denote the restriction of $W''$ to $V$. Then the argument above shows that the N\'{e}ron-Severi group of a general fibre of $W''_V \to V$ is fixed under the action of monodromy in $V$. Therefore, \cite[Proposition 4.1]{DHNT:2013} shows that we may undo the Kummer construction on $W''_V$ to obtain a threefold $\cA_V \to V$ fibred by products of elliptic curves $E_1 \times E_2$.

\subsection{Elliptic Surfaces} \label{subsection:elliptic surfaces}

The aim of the remainder of this section is to directly construct a birational model for $\cA_V \to V$ as a product of elliptic surfaces and to explain how this can be used to construct a birational model for $W \to \PP^1$.

Begin by letting $\cA_p$ be a fibre of $\cA_V \to V$ over a point $p \in V$, then $\cA_p$ is isomorphic to a product of smooth elliptic curves $E_1 \times E_2$. Furthermore, by Theorem \ref{theorem:M-polarized K3} and equations \eqref{eq:pi_Y} and \eqref{eq:sigma_Y}, the $j$-invariants of these elliptic curves are given by the roots of the quadratic equation
\begin{equation} j^2 - j + \frac{uv}{(u+v)^2 12^6 \xi_0} = 0, \label{eq:j quadratic} \end{equation}
where $[u:v] = f \circ g (p) \in \PP^1$.

We claim that $\mathcal{A}_V$ is isomorphic over $V$ to a fibre product of elliptic surfaces. This will be proved in several steps, the first of which is to show that the fibration on $\mathcal{A}_V$ admits a section.

\begin{lemma} \label{lemma:section of A} The fibration $\mathcal{A}_V \to V$ admits a section $\sigma\colon V \to \cA_V$. \end{lemma} 
\begin{proof} Note that the fibres of $W''_V$ over $V$ are all smooth Kummer surfaces. Furthermore, by the results of \ref{subsection:construction}, we see that the N\'{e}ron-Severi group of a fibre $W''_p$ of $W''_V$ is invariant under monodromy in $V$. The sixteen exceptional $(-2)$-curves arising from the Kummer construction on $W''_p$ are thus fixed under monodromy in $V$, so sweep out sixteen divisors on $W''_V$. On progression to $\cA_V$, these divisors are contracted to give sections of $\cA_V \to V$. \end{proof} 

Using this, we can now prove the following proposition.

\begin{proposition} \label{proposition:elliptic surfaces} $\mathcal{A}_V \to V$ is isomorphic over $V$ to a fibre product $\cE_1 \times_{C} \cE_2$ of minimal elliptic surfaces $\cE_{1,2} \to C$ with section. Furthermore, the $j$-invariants of the elliptic curves $E_1$ and $E_2$ forming the fibres of $\cE_1$ and $\cE_2$ over a point $p \in C$ are related by $j(E_1) + j(E_2) = 1$
\end{proposition}
\begin{proof} Let $\cA_p = E_1 \times E_2$ be a fibre of $\cA_V$ over a point $p \in V$. The section $\sigma\colon V \to \cA_V$ defined in Lemma \ref{lemma:section of A} intersects $\cA_p$ at a single point; we identify $E_1$ and $E_2$ with the two elliptic curves in $\cA_p$ passing through this point. Then we may deform $E_1$ and $E_2$ along paths in $\sigma(V)$ to obtain elliptic curves in the other fibres of $\cA_V$. Over any simply connected subset of $V$, this gives a decomposition of $\cA_V$ as a fibre product of two elliptic surfaces with section. To show that this decomposition extends to all of $V$, we have to show that the action of monodromy in $V$ cannot switch $E_1$ and $E_2$.

Let $\gamma$ be a path in $V$ that begins at $p \in \cA_p$, loops around one of the punctures in $V$, then returns to $p$. If we deform $E_1$ along $\sigma(\gamma)$, we obtain an elliptic curve $\gamma(E_1)$ in $\cA_p$ passing through $\sigma(p)$. Thus, $\gamma(E_1)$ must be either $E_1$ or $E_2$. 

We can determine which it is by tracking its $j$-invariant using equation \eqref{eq:j quadratic}. The discriminant of this quadratic vanishes singly at the points $[u:v] \in \{[\alpha:1],[\beta:1]\}$ and doubly at $[u:v] = [0:1]$, which are outside $V$, so we cannot have $j(E_1) = j(E_2)$ at $p$. We can therefore determine $\gamma(E_1)$ by comparing its $j$-invariant to $j(E_1)$ and $j(E_2)$.

If $\gamma$ is a loop around one of the punctures in $V$ above $\{[0:1],[-1:1]\} \subset \PP^1$, then the discriminant of \eqref{eq:j quadratic} does not vanish inside $f \circ g(\gamma)$ and we must have $j(\gamma(E_1)) = j(E_1)$, i.e. $\gamma(E_1) = E_1$. If $\gamma$ is a loop around one of the punctures in $V$ above $[1:0] \in \PP^1$, then $f \circ g(\gamma)$ loops three times around a point where the discriminant of \eqref{eq:j quadratic} vanishes doubly, so again we must have $j(\gamma(E_1)) = j(E_1)$ and $\gamma(E_1) = E_1$. Finally, if $\gamma$ is a loop around one of the punctures in $V$ above $\{[\alpha:1],[\beta:1]\} \subset \PP^1$, then $f \circ g(\gamma)$ loops twice around a point where the discriminant of \eqref{eq:j quadratic} vanishes singly, so we must have $j(\gamma(E_1)) = j(E_1)$ and $\gamma(E_1) = E_1$.

Thus, $\gamma(E_1) = E_1$ for any loop $\gamma$ in $V$. So $E_1$ sweeps out an elliptic surface $p_1\colon \cE_{1,V} \to V$ over $V$, with section given by $\sigma$. Similarly, $E_2$ also sweeps out an elliptic surface $p_2\colon \cE_{2,V} \to V$ with the same section, so we see that $\cA_V$ is isomorphic to the fibre product $\cE_{1,V} \times_{V} \cE_{2,V}$. 

Now, by \cite[Theorem 2.5]{Nakayama:1988}, there are unique extensions of $p_{i}\colon \cE_{i,V} \to V$ to minimal elliptic surfaces $p_i\colon \cE_i \to C$ over $C$, for $i = 1,2$. By construction, we have an isomorphism between $\cA_V$ and $\cE_1 \times_C \cE_2$ over $V$. Finally, the statement about the $j$-invariants is an easy consequence of equation \eqref{eq:j quadratic}.\end{proof}

Thus to construct a birational model for $\cA_V$, it is enough to construct the elliptic surfaces $\cE_1$ and $\cE_2$. In order to do this, we begin by studying their properties. Recall that, by the discussion in \ref{subsection:construction}, there is an action of $S_3 \times S_3$ on $C$. Let $i \colon C \to C$ be the involution defined by the transposition $(23)$. Then we have:

\begin{lemma} \label{lemma:elliptic isomorphism} The involution $i$ induces an isomorphism $\cE_1 \to \cE_2$.\end{lemma}
\begin{proof} As $\cE_1$ and $\cE_2$ are minimal, any birational map $\cE_1 \to \cE_2$ is an isomorphism by \cite[Proposition II.1.2]{Miranda:1989}. So suffices to show that $i$ induces an isomorphism over the open set $V$. To do this, we show that the fibre of $\cE_1$ over a point $p \in V$ is isomorphic to the fibre of $\cE_2$ over $i(p)$, for then we can define an isomorphism $\cE_{1,V} \to \cE_{2,V}$ by $(p,e) \mapsto (i(p), e)$, where $e$ is any point in the fibre of $\cE_1$ over $p$.

So let $E_k$ denote the fibre of $\cE_k$ over $p \in V$ and let $E_k'$ denote the fibre of $\cE_k$ over $i(p)$. By equation \eqref{eq:j quadratic}, we see that the $j$-invariants of $\{E_1,E_2\}$ are equal to those of $\{E_1',E_2'\}$. So to show that $E_1 \cong E_2'$,  we just have to show that $j(E_1) = j(E_2')$. 

To do this, we perform a calculation similar to the one used to prove Proposition \ref{proposition:elliptic surfaces}. Let $\gamma$ be a loop in $f\circ g(V)$ that begins at $f \circ g(p)$, loops once around $[\alpha:1] \in \PP^1$, then returns to $f \circ g(p)$. Then one of the preimages $\gamma'$ of $\gamma$ under $(f \circ g)$ is a path from $p$ to $i(p)$. As in the proof of Proposition \ref{proposition:elliptic surfaces}, if we deform $E_1,E_2 \subset \cA_p$ along $\sigma(\gamma')$ we obtain the fibres of $\cE_1,\cE_2$ over $i(p)$. Since $\gamma$ loops once around a point where the discriminant of \eqref{eq:j quadratic} vanishes singly, such a deformation swaps the $j$-invariants of the fibres of $\cE_1$ and $\cE_2$. So $j(E_1) = j(E_2')$ and $j(E_2) = j(E_1')$, as required.
\end{proof}

Using this lemma, we see that in order to construct a birational model for $\cA_V$ it is enough to construct a birational model for the elliptic surface $\cE_1$. As we can compute the fibrewise $j$-invariant of $\cE_1$ from equation \eqref{eq:j quadratic}, it just remains to study the forms of the singular fibres. These fibres will be calculated in the next section, but before that we conclude this section with a useful result about $\cE_1$.

\begin{lemma} \label{lemma:E1 automorphisms} Let $G$ denote the group of order $18$ obtained by intersecting $S_3 \times S_3 \subset S_6$ with the alternating group $A_6$. Then the action of $G$ on $C$ induces automorphisms of $\cE_1$. \end{lemma}
\begin{proof} The proof of this lemma proceeds using the same arguments used to prove Lemma \ref{lemma:elliptic isomorphism}.\end{proof} 

\subsection{Singular Fibres.} \label{subsection:degenerate}

The aim of this subsection is to compute the types of the singular fibres appearing in $\cE_1$. We begin by noting that, by equation \eqref{eq:j quadratic}, we have $j(\cE_1) = \infty$ with multiplicity two at the points lying over $[u:v] = [1:-1]$ and $j(\cE_1) \in \{0,1\}$ with multiplicity six at the points lying over $[u:v] \in \{[0:1],[1:0]\}$. Without loss of generality, we may assume that $j(\cE_1) = 1$ over $[u:v] = [0:1]$ and $j(\cE_1) = 0$ over $[u:v] = [1:0]$ (making the opposite choice simply corresponds to interchanging the labels of $\cE_1$ and $\cE_2$). This data determines the singular fibres of $\cE_1$ up to quadratic twists.

To determine these fibres exactly, we study the action of monodromy on the cohomology of the fibres of $\cE_1$. More precisely, let $p$ be a point in $V$ and let $E_i$ denote the fibre of $\cE_i$ over $p$. Let $\gamma$ be a loop in $V$ that starts at $p$, travels once around a puncture in $V$, then returns to $p$. Then monodromy around $\gamma$ acts on $H^1(E_1,\ZZ)$ and the precise form of this action determines the type of the singular fibre enclosed by $\gamma$.

Note that the fibre of $\cA_V \to V$ over $p$ is isomorphic to $E_1 \times E_2$. By the discussion at the end of \cite[Section 3.5]{ClingherDoran:2007}, we see that there is a canonical isomorphism
\[ H^1(E_1,\ZZ) \otimes H^1(E_2,\ZZ) \cong i(M)^{\perp} \subset H^2(Y_{1,p},\ZZ),\]
where $Y_{1,p}$ is the fibre of $Y_1 \to \PP^1$ over $f \circ g(p)$ (which is a smooth $M$-polarized K3 surface) and $i\colon M \to H^2(Y_{1,p},\ZZ)$ is the lattice embedding defining the $M$-polarisation on $Y_{1,p}$ (see Definition \ref{definition:M-polarized}).

Thus we can obtain information about the action of monodromy around $\gamma$ on $H^1(E_1,\ZZ)$ by studying the action of monodromy around $f \circ g(\gamma)$ on $i(M)^{\perp}$. To do this, we will study the Picard-Fuchs equation of $Y_1$. This equation has been computed in the degenerate case $\xi_0 = \frac{1}{12^6}$ by Chen, Doran, Kerr and Lewis in \cite[Section 5.2]{CDKL:2011}. In order to use their result in our computations, we perform a brief study of this case.

So suppose $\xi_0 = \frac{1}{12^6}$. The construction detailed in \ref{subsection:construction} proceeds much as before, except in this case we find that $\alpha = \beta = 1$. Due to this, the ramification points of the $36$-fold cover $g \colon C \to \PP^1$ lying over $[u:v] \in \{[\alpha:1],[\beta:1]\}$ collide and $g$ splits into two disjoint $18$-fold covers (which are preserved by even permutations in $S_3 \times S_3 \subset S_6$ and exchanged by odd ones). However, this splitting does not affect the ramification behaviour over $[u:v] \in \{[0:1],[-1:1],[1:0]\}$, so the singular fibres of $\cE_1$ and $\cE_2$ over these points are unaffected.

In this setting, Chen, Doran, Kerr and Lewis show that the Picard-Fuchs equation on $Y_1$ splits as a product of the second-order ODE's
\begin{align*}
f_1''(r) + \frac{3r+1}{2r(r+1)}f_1'(r) + \frac{5}{144r(r+1)}f_1(r) &= 0, \\
f_2''(r) + \frac{3r+1}{2r(r+1)}f_2'(r) + \frac{5}{144r^2(r+1)}f_2(r) &= 0,
\end{align*}
where $r = \frac{u}{v}$ is an affine parameter on the base $\PP^1$. Furthermore, this splitting corresponds exactly with the splitting of $i(M)^{\perp}$ as $H^1(E_1,\ZZ) \otimes H^1(E_2,\ZZ)$ so we can use it to study the singular fibres of $\cE_1$.

We recognize the first of these equations as the differential equation for the hypergeometric function $_2F_1(\frac{1}{12},\frac{5}{12};1\mid r+1)$. Matrix generators for the monodromy group with respect to this system are well-known (see, for instance, \cite[Section 2.4]{iksy:1991}), so we can immediately deduce the action of monodromy around a point $[u:v] \in \PP^1$ on $H^1(E_1,\ZZ)$. The results are shown in \ref{table:Y1monodromy}. 

\begin{table}
\begin{tabular}{|c|c|}
\hline 
Point $[u:v] \in \PP^1$ & Monodromy \\
\hline
$[0:1]$ & $i\left[ \begin{array}{cc} 0 & 1 \\ -1 & 0 \end{array} \right]$ \\
$[-1:1]$ & $\left[ \begin{array}{cc} 1 & 1 \\ 0 & 1 \end{array} \right]$ \\
$[1:0]$ & $i\left[ \begin{array}{cc} 0 & 1 \\ -1 & -1 \end{array} \right]$ \\
\hline
\end{tabular}
\caption{Action of monodromy around a point $[u:v] \in \PP^1$ on $H^1(E_1,\ZZ)$.}
\label{table:Y1monodromy}
\end{table}

Note, however, that this table does not give the action of monodromy in $C$ on $H^1(E_1,\ZZ)$. To obtain that, we have to take account of the cover $(f \circ g)\colon C \to \PP^1$. Let $\gamma$ be a loop in $C$ around one of the preimages $(f \circ g)^{-1}[0:1]$. Then $f \circ g(\gamma)$ is a loop in $\PP^1$ that encircles $[0:1]$ six times, so the action of monodromy around $\gamma$ on $H^1(E_1,\ZZ)$ is as the action of a loop encircling $[0:1] \in \PP^1$ raised to the sixth power. Similarly, the action of monodromy around a loop in $C$ about one of the preimages $((f \circ g)^{-1}[-1:1]$ (resp. $((f \circ g)^{-1}[1:0]$) is as the action of a loop encircling $[-1:1] \in \PP^1$ (resp. $[1:0] \in \PP^1$) raised to the power of two (resp. six). The resultant monodromy matrices and the corresponding singular fibres of $\cE_1$ are given in \ref{table:E1monodromy}.

\begin{table}
\begin{tabular}{|c|c|c|c|c|}
\hline 
$f \circ g(p)$ & $j$-invariant & Multiplicity of $j$ & Monodromy & Fibre \\
\hline
$[0:1]$ & $1$ & $6$ & $\left[ \begin{array}{cc} 1 & 0 \\ 0 & 1 \end{array} \right]$ & $I_0$ \\
$[-1:1]$ & $\infty$ & $2$ & $\left[ \begin{array}{cc} 1 & 2 \\ 0 & 1 \end{array} \right]$ & $I_2$ \\
$[1:0]$ & $0$ & $6$ & $\left[ \begin{array}{cc} -1 & 0 \\ 0 & -1 \end{array} \right]$ & $I_0^*$ \\
\hline
\end{tabular}
\caption{Singular fibres of $\cE_1$ over points $p \in C$.}
\label{table:E1monodromy}
\end{table}

\subsection{The forward construction.} \label{subsection:reconstructing}

The last thing we will discuss is how to use this information to construct a birational model for $W \to \PP^1$ and what this process can tell us about $W$.

Our beginning data is the two covers $C \stackrel{g}{\to} \PP^1 \stackrel{f}{\to} \PP^1$ (which can be determined from $\xi_0$) and the elliptic surface with section $\cE_1$ on $C$. Note that $S_3 \times S_3$ acts on $C$ and that the group $G$ of order $18$ given by the intersection of $S_3 \times S_3$ with $A_6$ acts on $\cE_1$.

We first construct $\cE_2$ as $i(\cE_1)$, where $i\colon C \to C$ is the involution defined by the transposition $(23)$ in $S_3 \times S_3$. Then we perform the fibrewise Kummer construction on the fibre product $\cE_1 \times_C \cE_2$ to obtain a birational model $\cF$ for $W''$.

The next step is to perform a quotient by $S_3 \times S_3$ to obtain a birational model for $W'$. However, the action of $S_3 \times S_3$ is not the obvious one induced by the action of $S_3 \times S_3$ on $\cE_1 \times_C \cE_2$ (if it were, we would be able to undo the Kummer construction on $W'$, which we have previously shown to be impossible). Instead, we compose this action with the fibrewise automorphism induced by the action of $S_3 \times S_3$ on the sets of curves $\{F_3,F_4,F_5\}$ and $\{F_6,F_7,F_8\}$.

\begin{remark} As the fibrewise Kummer construction defines a natural double Kummer pencil on smooth fibres of $\cF$, we can use the results of Kuwata and Shioda \cite[Section 5.2]{KuwataShioda:2008} to define the elliptic fibration $\Psi$ on $\cF$. The curves $\{F_3,\ldots,F_8\}$ are then the components of the $I_2$ fibres that do not meet a chosen section.

Once $\{F_3,\ldots,F_8\}$ are known, the automorphisms permuting them may be calculated explicitly as compositions of the symplectic automorphisms $f_{r'}$ from Keum and K\={o}ndo \cite[Section 4.1]{keumkondo:2001}. We will give an example of such calculations in Example \ref{example:ExE fibres}. \end{remark}

Quotienting by this $S_3 \times S_3$ action, we obtain a birational model for $W' \to \PP^1$. Finally, there is a natural involution identifying the fibres over $[s:t]$ and $[-s:t]$ in this model. Quotienting by this involution, we obtain the model for $W \to \PP^1$ we desire.

We conclude this section with an example showing how the action of $S_3 \times S_3$ can be computed around a singular fibre; this will allow us to gather information about some of the singular fibres in $W$. For simplicity we consider the fibres over $[\alpha:1]$ and $[\beta:1]$; the remaining fibres have too many components to allow such simple explicit calculations to be performed.

\begin{example} \label{example:ExE fibres} We will calculate the form of the singular fibre in our birational model over $[\alpha:1]$ (the calculation for $[\beta:1]$ is identical). Let $\Delta$ denote a small disc in $\PP^1$ centred at $\alpha$, and let $\Delta'$ denote one of the connected components of $(f \circ g)^{-1}(\Delta)$. Then the map $(f \circ g)\colon \Delta' \to \Delta$ is a double cover, ramified over $0 \in \Delta$.

Over $\Delta'$ the elliptic surface $\cE_1$ has no singular fibres. The involution $i$ acts on $\Delta'$ to exchange the sheets of the double cover $(f \circ g)$. Define $\cE_2 := i(\cE_1)$. Note that the fibres of $E_1$ and $E_2$ over $(f \circ g)^{-1}(0) \in \Delta'$ are isomorphic.

Perform the fibrewise Kummer construction on $\cE_1 \times_{\Delta'} \cE_2$ to obtain $\cF \to \Delta'$. Each fibre of $\cF$ comes equipped with $24$ distinguished $(-2)$-curves forming a \keyterm{double Kummer pencil} $\{G_i,H_j,E_{ij} \mid 0 \leq i,j \leq 3\}$ (as defined in \cite[Definition 3.18]{ClingherDoran:2007}). As the involution $i$ exchanges $\cE_1$ and $\cE_2$, it acts on this double Kummer pencil to exchange $G_i \leftrightarrow H_i$ and $E_{ij} \leftrightarrow E_{ji}$.

Now use \cite[Section 5.2]{KuwataShioda:2008} to define the elliptic fibration $\Psi$ on $\cF$. In this fibration, the curves $\{F_3,\ldots,F_8\}$ are the components of the $I_2$ fibres that do not meet a chosen section (which we take to be $G_2$). This set divides naturally into two triples, $\{F_3,F_4,F_5\}$ and $\{F_6,F_7,F_8\}$. By \cite[Section 5.2]{KuwataShioda:2008}, the curves $\{F_3,F_4,F_5\}$ (up to labelling) are given in the double Kummer pencil by
\begin{align*}
F_3 &\sim E_{33} \\
F_4 &\sim \frac{1}{2}\left(\sum_{i=0}^3 \left(G_i + H_i\right) + \sum_{i,j} E_{ij} - 2E_{00} - 2E_{12} - 2E_{21}\right)\\
F_5 &\sim  \sum_{i=0}^3 \left(G_i + H_i\right) + \sum_{i,j} E_{ij} - 2E_{00} - E_{13} - E_{31} - E_{12} - E_{21} - E_{22}
\end{align*}
For simplicity, we will assume that the fixed curves swapped by monodromy around $0\in \Delta$ are $F_3$ and $F_4$ (the other choices give the same result, but the calculations are substantially longer). By the results of \cite[Section 4.1]{keumkondo:2001}, this switching is realized by a symplectic involution $\varphi$ of the Kummer surface that takes
\begin{align*}
\varphi(G_i) &= H_i \\
\varphi(H_i) &= G_i \\
\varphi(E_{ij}) &= \left\{\begin{array}{ll} E_{ji} & \mathrm{if}\ (i,j) \notin \{(0,0),(1,2),(2,1),(3,3)\} \\
 E_{ji} + D & \mathrm{if}\ (i,j) \in \{(0,0),(1,2),(2,1),(3,3)\} \end{array} \right.
\end{align*} 
where $D$ denotes the divisor
\[ D:= \frac{1}{2}\left(\sum_{i=0}^3 \left(G_i + H_i\right) + \sum_{i,j} E_{ij} - 2E_{00} - 2E_{12} - 2E_{21} - 2E_{33}\right).\]

The involution defining the map from $\cF$ to the birational model for $W$ is realized as the composition of $i$ with the fibrewise involution given by $\varphi$. The fixed points of this involution all lie in the fibre of $\cF$ over $(f \circ g)^{-1}(0)$. This fibre is isomorphic to $\mathrm{Kum}(E \times E)$ for some smooth elliptic curve $E$.

In the fibre of $\cF$ over $(f \circ g)^{-1}(0)$, the $I_2$ fibres in $\Psi$ that are exchanged by the fibrewise involution $\varphi$ collide to give an $I_4$ fibre. Let $D_1$ and $D_2$ denote the two disjoint $(-2)$-curves in this $I_4$ fibre that do not meet the sections $\{G_2,G_3,H_2,H_3\}$. Then $D = D_1 + D_2$, but $D_1$ and $D_2$ cannot be expressed as a combination of divisors from the double Kummer pencil $\{G_i,H_j,E_{ij}\}$ (the N\'{e}ron-Severi group of a generic fibre of $\cF$ has rank $18$ and is spanned by the divisors in the double Kummer pencil, but this rank jumps to $19$ on the special fibre over $(f \circ g)^{-1}(0)$ and the double Kummer pencil no longer spans). The involution $\varphi$, which is defined in terms of the double Kummer pencil, is not well-defined on $D_1$ and $D_2$, so they must be contracted before we can quotient.

After performing this contraction, the threefold $\cF$ has two nodes and its fibre over $(f \circ g)^{-1}(0)$ has two $A_1$ singularities. The involution $\varphi$ acts on the remaining divisors in $\mathrm{Kum}(E \times E)$ to exchange $G_i \leftrightarrow H_i$ and $E_{ij} \leftrightarrow E_{ji}$. But this is precisely the same as the action of the involution $i$. So the composition $(\varphi \circ i)$ acts trivially on the fibre of $\cF$ over $(f \circ g)^{-1}(0)$.

After performing the quotient by $(\varphi \circ i)$, we find that the threefold total space is smooth over the disc $\Delta$. Its fibre over $0 \in \Delta$ is isomorphic to the singular K3 surface obtained from $\mathrm{Kum}(E \times E)$ by contracting the two $(-2)$-curves $D_1$ and $D_2$. This fibre therefore has two $A_1$ singularities.
\end{example}

\section{Moduli Spaces and the Mirror Map}\label{section:mirrormoduli}

The aim of this final section is to discuss the mirrors to the families $Y$ and $Z$, along with their subfamilies $Y_1$, $Z_2$ and $Z_1$. As has already been discussed in \ref{subsection:CI model} and \ref{subsection:AH model}, the mirror $Y^{\circ}$ of $Y$ is known to be a family of $(2,12)$-complete intersections inside a blow-up of $\WP(1,1,1,1,4,6)$ and the mirror $Z^{\circ}$ of $Z$ is known to be a family of degree $24$ hypersurfaces in a blow-up of $\WP(1,1,2,8,12)$.

We will exhibit a candidate mirror for the degeneration of $Z$ to $Z_2$ and $Z_1$ in terms of the K\"{a}hler moduli of $Z^{\circ}$, and discuss the mirror of the geometric transition between $Y$ and $Z$. This will enable us to match the components of the discriminant locus in the complex moduli space of $Z$ to the boundary components in the K\"{a}hler cone of $Z^{\circ}$. These considerations finally lead us to conclude that $Y_1$ and its mirror provide a counterexample to a conjecture of Morrison  \cite{Morrison:1999}.

\subsection{Complex and K\"{a}hler Moduli}

We begin by discussing the complex moduli of $Z$ and the K\"{a}hler moduli of $Z^{\circ}$. As noted in \ref{subsection:AH model}, the complex moduli are controlled by four parameters, $B$, $\psi_0$, $\psi_1$ and $\psi_s$. The discriminant locus was originally calculated by \cite{Billo:1998} and is given in \ref{subsection:singular hypersurfaces}. It splits into three parts: 
\begin{itemize}
\item The locus $B=0$, which we call $S_0$. From Equations \eqref{eq:pi_Z} and \eqref{eq:sigma_Z}, we see that this corresponds to the locus where the $M$-polarized K3 fibration on $Z$ becomes isotrivial.
\item The locus $\psi_s = \pm B$, which we call $S_b$. This corresponds to the degenerate family $Z_2$.
\item The loci defined by the relations $(\psi_0^6 + \psi_1)^2 + \psi_s = \pm B$ and $\psi_1^2 + \psi_s = \pm B$. These two relations are switched by an appropriate change of coordinates, so only give a single locus in the complex moduli space, which we call $S_a$.
\end{itemize}
Note that the intersection $S_a \cap S_b$ corresponds to the degenerate family $Z_1$.

Next we discuss the K\"{a}hler moduli of $Z^{\circ}$. This has been studied by both Scheidegger \cite[Appendix C.3]{Scheidegger:2001} and Hosono, Klemm, Theisen and Yau \cite{HKTY:1997}; we will use Scheidegger's notation in the following discussion. 

A generic degree $24$ hypersurface in $\WP(1,1,2,8,12)$ has an elliptic curve $C$ of $A_1$ singularities along with an exceptional $\ZZ_4$ point lying on this curve. The mirror $Z^{\circ}$ of $Z$ is given by the crepant resolution of these singularities and contains two exceptional components: a ruled surface $E$ over the curve $C$ and a Hirzebruch surface $F \cong \mathbb{F}_2$ coming from the blow up of the exceptional point. Along with the hyperplane section $H$, these classes span the K\"{a}hler cone of $Z^{\circ}$.

By \cite[Equation (A.40)]{HKTY:1997}, the K\"{a}hler cone of $Z^{\circ}$ is given in terms of these classes by
\[ \{t_H H + t_E E + t_F F \mid t_H + t_F > 0,\ t_E < 0,\ t_E - 2t_F > 0\}.\]
Scheidegger shows that a basis of the Mori cone of $Z^{\circ}$ is given by the classes $h$, $d$ and $l$, where $h$ is the class of a section of $E$ (and is dual to $H$), $d$ is the class of a fiber of $F$, and $l$ is the class of a fiber of $E$ (which is also the class of the $(-2)$-section of $F$). Proceeding to a boundary component of the K\"{a}hler cone corresponds to contracting one of these classes, as follows:
\begin{itemize}
\item The boundary component $t_E = 0$ is spanned by $H$ and $F$. Scheidegger shows that $H.l = F.l = 0$, so $l$ is contracted on this component. This corresponds to the birational morphism contracting $E$ along its ruling onto the curve $C$. The $(-2)$-section in $F$ is also contracted, giving a cone $\mathbb{F}_2^0$.
\item The boundary component $t_E - 2t_F = 0$ is spanned by $H$ and $F + 2E$. Scheidegger shows that $H.d = (F + 2E).d = 0$, so $d$ is contracted on this component. This corresponds to the birational morphism contracting $F$ along its ruling.
\item The boundary component $t_H + t_F = 0$ is spanned by $E$ and $H-F$. Scheidegger shows that $E.h = (H-F).h = 0$, so $h$ is contracted on this component. This does not give rise to a birational contraction.
\end{itemize}

\subsection{The Mirror Map} \label{section:mirror}

Let $Z^{\circ}_2$ denote the threefold obtained by contracting the divisor $E$ in $Z^{\circ}$ along its ruling; it may also be seen as the threefold obtained by blowing-up the $\ZZ_4$ point in a generic degree $24$ hypersurface in $\WP(1,1,2,8,12)$. Along with the hyperplane class $H$, it contains one exceptional component $F$, which is isomorphic to a cone $\mathbb{F}_2^0$. As explained above, the K\"{a}hler modulus of $Z^{\circ}_2$ lies on the boundary component $t_E = 0$, which is spanned by the classes of $H$ and $F$.

Now let $Z^{\circ}_1$ denote the threefold obtained from $Z^{\circ}_2$ by contracting $F$. Then $Z^{\circ}_1$ is isomorphic to a generic degree $24$ hypersurface in $\WP(1,1,2,8,12)$. Its K\"{a}hler modulus lies on the ray given by the intersection of the boundary components $t_E = 0$ and $t_E - 2t_F = 0$, which is spanned by the class of $H$, the hyperplane section.

We claim that $Z^{\circ}_2$ and $Z^{\circ}_1$ are mirror to $Z_2$ and $Z_1$ respectively. To justify this, we exhibit a mirror for the geometric transition between $Y$ and $Z$.

Recall from \ref{subsection:desingularization} that we may obtain a geometric transition from $Z$ to $Y$ by first degenerating to $Z_2$, then blowing-up the curve of singularities in $Z_2$. Furthermore, this transition extends to $Z_1$: we may degenerate $Z_2$ further to $Z_1$, then blow-up the curve of singularities in $Z_1$ to obtain the (singular) threefold $Y_1$.  This relationship was summarized by Diagram \eqref{equation:mirror1}.

Next, we discuss the mirror picture. Note first that we may embed $Z^{\circ}_1$ into $\WP(1,1,1,1,4,6)$ as a (non-generic) complete intersection of type $(2,12)$, via the degree two Veronese embedding $\WP(1,1,2,8,12) \hookrightarrow \WP(2,2,2,2,8,12)$. This admits a partial smoothing to a generic complete intersection of type $(2,12)$ in $\WP(1,1,1,1,4,6)$, which should be thought of as mirror to $Y_1$ and will be denoted by $Y_1^{\circ}$; note that this exhibits $Y_1$  as mirror to a complete intersection of type $(2,12)$ in $\WP(1,1,1,1,4,6)$, precisely as expected. The threefold $Y_1^{\circ}$ has an isolated Gorenstein canonical singularity, which may be blown-up once to obtain the smooth Calabi-Yau threefold $Y^{\circ}$, which is the mirror of $Y$. Kobayashi \cite[Theorem 1]{Kobayashi:1998} has shown that the exceptional locus of this resolution is isomorphic to $\PP^1 \times \PP^1$. 

By considering the blow-up of a family of threefolds in $\WP(1,1,1,1,4,6)$, we may naturally see $Z_2^{\circ}$ as a degeneration of $Y^{\circ}$, where the generic $\PP^1 \times \PP^1$ exceptional locus degenerates to $\mathbb{F}_2^0$ in the limit. We thus have a diagram, where wiggly arrows denote degeneration and straight arrows denote blow-ups:
\begin{equation}\label{equation:mirror2}\xymatrix{ & Y^{\circ} \ar[r] \ar@{~>}[d] & Y_1^{\circ} \ar@{~>}[d] \\
Z^{\circ} \ar[r] & Z_2^{\circ} \ar[r] & Z_1^{\circ}
}\end{equation}
Diagrams \eqref{equation:mirror1} and \eqref{equation:mirror2} should be thought of as mirror to one another. In particular, we see that there is a mirror geometric transition from $Z^{\circ}$ to $Y^{\circ}$, given by first contracting to $Z_2^{\circ}$, then smoothing to $Y^{\circ}$. In fact, this mirror correspondence between the geometric transitions from $Y$ to $Z$ and from $Z^{\circ}$ to $Y^{\circ}$ provides an explicit example of the kind of mirror correspondence discussed by Mavlyutov \cite[Section 5]{Mavlyutov:2011}.

From this, we can deduce the mirror correspondence between the components of the discriminant locus in the complex moduli space of $Z$ and the components of the boundary of the K\"{a}hler cone of $Z^{\circ}$. We find:
\begin{itemize}
\item The locus $S_b$ in the complex moduli space of $Z$, which corresponds to the degenerate family $Z_2$, must be mirror to the boundary component $t_E = 0$ in the K\"{a}hler moduli space of $Z^{\circ}$, which corresponds to $Z_2^{\circ}$.
\item The locus $S_a$ in the complex moduli space of $Z$ must be mirror to the boundary component $t_E - 2t_F = 0$ in the K\"{a}hler moduli space of $Z^{\circ}$, as the intersection $S_a \cap S_b$ corresponds to the degenerate family $Z_1$, and the mirror family $Z_1^{\circ}$ corresponds to the boundary ray $t_E = t_E - 2t_F = 0$.
\item By elimination, the remaining locus $S_0$  in the complex moduli space of $Z$ must be mirror to the boundary component $t_H + t_F = 0$ in the K\"{a}hler moduli space of $Z^{\circ}$.
\end{itemize}

We conclude with a brief note on a conjecture of Morrison \cite{Morrison:1999}, which states that a singular Calabi-Yau space has a Calabi-Yau resolution if and only if its mirror has a Calabi-Yau smoothing. However, in our case we find that, by a theorem of Gross \cite[Theorem 5.8]{Gross:1997}, the fact that the exceptional locus in the mirror $Y^{\circ}$ is isomorphic to $\PP^1 \times \PP^1$ implies that $Y^{\circ}_1$ has a Calabi-Yau smoothing. But, by Corollary \ref{corollary:Y1sympres}, $Y_1$ does not admit a Calabi-Yau resolution. This provides a counterexample to Morrison's conjecture.


\bibliography{Bibliography}
\bibliographystyle{alpha}

\end{document}